\documentclass[oneside,english,reqno]{amsart}
\usepackage{mathptmx}
\usepackage[T1]{fontenc}
\usepackage[latin9]{inputenc}
\usepackage{babel}

\usepackage{textcomp}
\usepackage{amsthm}
\usepackage{graphicx}
\usepackage{amssymb}
\usepackage[unicode=true, pdfusetitle,
 bookmarks=true,bookmarksnumbered=false,bookmarksopen=false,
 breaklinks=false,pdfborder={0 0 1},backref=false,colorlinks=false]
 {hyperref}

\makeatletter

\providecommand{\tabularnewline}{\\}

\theoremstyle{plain}
\newtheorem{thm}{Theorem}[section]
  \theoremstyle{plain}
  
  \theoremstyle{definition}
  \newtheorem{defn}[thm]{Definition}
  \theoremstyle{plain}
  \newtheorem{lem}[thm]{Lemma}
  \theoremstyle{plain}
  \newtheorem{fact}[thm]{Fact}
 \theoremstyle{definition}
  \newtheorem{example}[thm]{Example}
  \theoremstyle{remark}
  \newtheorem{rem}[thm]{Remark}
  \theoremstyle{plain}
  \newtheorem{conjecture}[thm]{Conjecture}

\newcommand{\lip}{\mbox{\rm lip}\,}
\newcommand{\dom}{\mbox{\rm dom}}
\newcommand{\epi}{\mbox{\rm epi}\,}
\newcommand{\gph}{\mbox{\rm Graph}}
\newcommand{\cl}{\mbox{\rm cl}}
\newcommand{\co}{\mbox{\rm co}}

\numberwithin{equation}{section}
\numberwithin{figure}{section}
\title[Subdifferential analysis of differential inclusions via discretization]{Subdifferential analysis of differential inclusions via discretization}

\makeatother

\begin{document}

\author{C.H. Jeffrey Pang}

\curraddr{Massachusetts Institute of Technology, Department of Mathematics,
2-334, 77 Massachusetts Avenue, Cambridge MA 02139-4307.}

\email{chj2pang@mit.edu}



\date{\today{}}
\begin{abstract}
 The framework of differential inclusions encompasses modern optimal
control and the calculus of variations. Necessary optimality conditions
in the literature identify potentially optimal paths, but do not show
how to perturb paths to optimality. We first look at the corresponding
discretized inclusions, estimating the subdifferential dependence
of the optimal value in terms of the endpoints of the feasible paths.
Our approach is to first estimate the coderivative of the reachable
map. The discretized (nonsmooth) Euler-Lagrange and transversality
conditions follow as a corollary. We obtain corresponding results
for differential inclusions by passing discretized inclusions to the
limit.
\end{abstract}
\maketitle
\tableofcontents{}

\section{Introduction}

The subject of this paper is the analysis of discretized differential
inclusions by calculating the coderivatives of the discretized reachable
map. We then pass these results to the limit to obtain results on
differential inclusions. We say that $S$ is a \emph{set-valued map
}or a \emph{multifunction}, denoted by $S:X\rightrightarrows Y$,
if $S(x)\subset Y$ for all $x\in X$. For $F:[0,T]\times\mathbb{R}^{n}\rightrightarrows\mathbb{R}^{n}$
and $C\subset\mathbb{R}^{n}\times\mathbb{R}^{n}$, consider the \emph{differential
inclusion}:

\begin{eqnarray}
 & \underset{x(\cdot)\in AC([0,T],\mathbb{R}^{n})}{\min} & \varphi\big(x(0),x(T)\big)\label{eq:main-diff-in}\\
 & \mbox{s.t. } & x^{\prime}(t)\in F\big(t,x(t)\big)\mbox{ for }t\in[0,T]\mbox{ a.e.}\nonumber \end{eqnarray}
Here, $AC([0,T],\mathbb{R}^{n})$ is the set of absolutely continuous
functions $x:[0,T]\to\mathbb{R}^{n}$. The constraint \[
(x(0),x(T))\in C\subset\mathbb{R}^{n}\times\mathbb{R}^{n}\]
is sometimes included, but this constraint can be easily incorporated
into the objective function $\varphi$. More details on differential
inclusions can be obtained in the texts \cite{AC84,AF90,Cla83,Mor06,Smirnov02,Vinter00}.
As is popularized in these texts, the differential inclusion framework
\eqref{eq:main-diff-in} encompasses optimal control and the calculus
of variations. 

In order to optimize \eqref{eq:main-diff-in}, much attention has
focused on necessary optimality conditions for a path $x(\cdot)$.
Such research were undertaken in the last few decades by Clarke, Loewen,
Rockafellar, Ioffe, Vinter, Mordukhovich, Kaskosz and Lojasiewicz,
Milyutin, Smirnov, Zheng, Zhu and others, building on results in the
calculus of variations and optimal control. For a history of the development
of the necessary optimality conditions, we refer to the previously
mentioned texts. The following conditions are currently understood
as useful necessary optimality conditions for a feasible path $\bar{x}(\cdot)$
of \eqref{eq:main-diff-in}:
\begin{itemize}
\item [(TC)](Transversality Condition) \\
$\qquad\qquad$$(-p(0),p(T))\in\partial\varphi(\bar{x}(0),\bar{x}(T))$
\item [(EL)](Euler-Lagrange Condition) \\
$\qquad\qquad$$p^{\prime}(t)\in-\overline{\co}D_{x}^{*}F\big(t,\bar{x}(t)\mid\bar{x}^{\prime}(t)\big)\big(p(t)\big)$
for $t\in[0,T]$ a.e.. 
\item [(WP)](Weierstrass-Pontryagin Maximum Principle) \\
$\qquad\qquad$$\left\langle -p(t),v-\bar{x}^{\prime}(t)\right\rangle \leq0$
for all $v\in F\big(t,\bar{x}(t)\big)$, $t\in[0,T]$ a.e..
\end{itemize}
While such necessary conditions are helpful in finding candidates
for a minimizing path, the deficiency in such necessary conditions
is that they do not give an indication on how to perturb a feasible
path to optimality. As a first step, we study the discrete inclusions
corresponding to the differential inclusion and calculate the dependence
of the differential inclusion on its initial point. 

Define the \emph{reachable map} (or \emph{attainable map}) $R:\mathbb{R}^{n}\rightrightarrows\mathbb{R}^{n}$
by \begin{eqnarray}
R(x_{0}) & := & \{y:\exists x(\cdot)\in AC([0,T],\mathbb{R}^{n})\mbox{ s.t. }\label{eq:R-diff-reachable}\\
 &  & \qquad x^{\prime}(t)\in F\big(t,x(t)\big)\mbox{ for }t\in[0,T]\mbox{ a.e.},\nonumber \\
 &  & \qquad x(0)=x_{0}\mbox{ and }x(T)=y\}.\nonumber \end{eqnarray}
In order to study \eqref{eq:main-diff-in}, we study \begin{eqnarray}
f(x):= & \underset{x}{\min} & \varphi(x,y)\label{eq:reachable-diff-in}\\
 & \mbox{s.t.} & y\in R(x)\nonumber \end{eqnarray}
We study \eqref{eq:reachable-diff-in} under the broader framework
of marginal functions. For a set-valued map $G:X\rightrightarrows Y$
and a function $\varphi:X\times Y\to\mathbb{R}$, the \emph{marginal
function }$f:X\to\mathbb{R}$ is \begin{equation}
f(x):=\inf\{\varphi(x,y):y\in G(x)\}.\label{eq:marginal-1}\end{equation}
One can view the value $x$ as a parameter of an optimization problem
in terms of $y$. A well studied example of a set-valued map $G$
is the map $G:\mathbb{R}^{n}\rightrightarrows\mathbb{R}^{m}$ defined
by \begin{eqnarray*}
G(x) & = & \{y\mid y\in F(x)+[\{0\}\times\mathbb{R}_{-}^{m_{2}}]\},\\
 &  & \qquad\mbox{where }m_{1}+m_{2}=m\mbox{ and }F:\mathbb{R}^{n}\to\mathbb{R}^{m}\mbox{ is smooth}\\
 & = & \{y\mid y_{i}=F_{i}(x)\mbox{ for }1\leq i\leq m_{1}\mbox{ and }y_{i}\leq F_{i}(x)\mbox{ for }m_{1}+1\leq i\leq m\}.\end{eqnarray*}
The sensitivity analysis of marginal functions can be analyzed with
tools of set-valued analysis. We denote the composition $S_{2}\circ S_{1}:X\rightrightarrows Z$
of set-valued maps $S_{1}:X\rightrightarrows Y$ and $S_{2}:Y\rightrightarrows Z$
in the usual way by \[
S_{2}\circ S_{1}(x)=\bigcup_{y\in S_{1}(x)}S_{2}(y).\]
Denote the epigraphical mapping of $\varphi$ and $f$ by $E_{\varphi}:X\times Y\rightrightarrows\mathbb{R}$
and $E_{f}:X\rightrightarrows\mathbb{R}$ respectively. Then $E_{\varphi}$
and $E_{f}$ satisfy the relation \begin{equation}
E_{f}(x)=E_{\varphi}\circ\bar{G}(x),\label{eq:epi-compose}\end{equation}
where $\bar{G}:X\rightrightarrows X\times Y$ is defined by $\bar{G}(x)=\{x\}\times G(x)$.
The relationship \eqref{eq:epi-compose} and a set-valued chain rule
can be used to express differentiability properties of $f$ in terms
of the coderivatives of $G$ and $\varphi$.

\subsection{Contributions of this paper}

In this work, we focus on the subdifferential analysis of the discretized
differential inclusion problem by finding $\partial f(x)$, where
$f$ is the discretized analogue of \eqref{eq:reachable-diff-in}.
Our approach is to look at the marginal function framework and calculate
the coderivatives of the reachable map $R(\cdot)$. The coderivative
of the reachable map gives new insight on the Euler-Lagrange Condition
(EL). We also study the limitations of a discrete analogue of the
Weierstrass-Pontryagin Maximum Principle (WP).

For a set valued map $S:\mathbb{R}^{n}\rightrightarrows\mathbb{R}^{m}$
between finite dimensional spaces, \cite{chardiff} recently established
that the convexified limiting coderivative characterizes the set of
positively homogeneous maps that are generalized derivatives of $S$
as defined in \cite{set_diff}. We will recall on this relation in
Section \ref{sec:prelim}, limiting our analysis to the finite dimensional
case. By making use of this result, we can obtain the convexified
limiting coderivative of $R(\cdot)$ by passing a sequence of discrete
problems to the limit. The marginal function framework allows us to
calculate the subdifferential dependence of the differential inclusion
in terms of its initial value.

\subsection{Outline}

In Section \ref{sec:prelim}, we recall standard definitions in variational
analysis and some results in \cite{chardiff} that will be used in
the later part of the paper. In Section \ref{sec:calculus-co-coderv},
we recall chain rules for coderivatives, and show how these results
can be easily extended for the convexified limiting coderivative.
In Section \ref{sec:discrete-incl}, we study the discretized differential
inclusion problem. Finally, in Section \ref{sec:cts-incl}, we study
the continuous inclusion problem by passing the discretized problems
in Section \ref{sec:discrete-incl} to the limit, and find formulas
for the convexified limiting coderivative of the reachable map.

\section{\label{sec:prelim}Preliminaries and notation}

This section recalls some standard definitions in variational analysis
and some other results in \cite{chardiff} that will be used in the
remainder of this paper. The texts \cite{RW98,Mor06} contain many
standard definitions in variational analysis, like inner and outer
semicontinuity (isc and osc) and the Pompieu-Hausdorff distance $\mathbf{d}(\cdot,\cdot)$.
We highlight some of definitions used most often in this paper. We
denote the set $\{1,2,\dots,N\}$ by $\overline{1,N}$. For set-valued
maps $H_{i}:\mathbb{R}^{n}\rightrightarrows\mathbb{R}^{m}$, $i=1,2$,
we let $H_{1}\subset H_{2}$ denote $H_{1}(x)\subset H_{2}(x)$ for
all $x$, or equivalently $\gph(H_{1})\subset\gph(H_{2})$. 

We recall the definition of coderivatives. 
\begin{defn}
\label{def:normal-cones}(Normal cones) For a set $C\subset\mathbb{R}^{n}$,
the \emph{regular normal cone} at $\bar{x}$ is defined as \[
\hat{N}_{C}(\bar{x}):=\{y\mid\left\langle y,x-\bar{x}\right\rangle \leq o(\|x-\bar{x}\|)\mbox{ for all }x\in C\}.\]
The \emph{limiting }(or\emph{ Mordukhovich})\emph{ normal cone} $N_{C}(\bar{x})$
is defined as $\limsup_{x\xrightarrow[C]{}\bar{x}}\hat{N}_{C}(x)$,
or as \[
N_{C}(\bar{x})=\{y\mid\mbox{there exists }x_{i}\xrightarrow[C]{}\bar{x},\, y_{i}\in\hat{N}_{C}(x_{i})\mbox{ such that }y_{i}\to y\}.\]

\begin{defn}
(Coderivatives) For a set-valued map $S:\mathbb{R}^{n}\rightrightarrows\mathbb{R}^{m}$
locally closed at $(\bar{x},\bar{y})\in\gph(S)$, the \emph{regular
coderivative} at $(\bar{x},\bar{y})$, denoted by $\hat{D}^{*}S(\bar{x}\mid\bar{y}):\mathbb{R}^{m}\rightrightarrows\mathbb{R}^{n}$,
is defined by\begin{align*}
v\in\hat{D}^{*}S(\bar{x}\mid\bar{y})(u) & \Leftrightarrow(v,-u)\in\hat{N}_{\scriptsize\gph(S)}(\bar{x},\bar{y})\\
 & \Leftrightarrow\left\langle (v,-u),(x,y)-(\bar{x},\bar{y})\right\rangle \leq o\big(\|(x,y)-(\bar{x},\bar{y})\|\big)\\
 & \qquad\qquad\mbox{ for all }(x,y)\in\gph(S).\end{align*}
The \emph{limiting }(or \emph{Mordukhovich) coderivative} at $(\bar{x},\bar{y})\in\gph(S)$
is denoted by $D^{*}S(\bar{x}\mid\bar{y}):\mathbb{R}^{m}\rightrightarrows\mathbb{R}^{n}$
and is defined by  \[
v\in D^{*}S(\bar{x}\mid\bar{y})(u)\Leftrightarrow(v,-u)\in N_{\scriptsize\gph(S)}(\bar{x},\bar{y}).\]
The \emph{convexified limiting coderivative} $\overline{\co}D^{*}S(\bar{x}\mid\bar{y}):\mathbb{R}^{m}\rightrightarrows\mathbb{R}^{n}$
is defined in the natural manner.
\end{defn}
\end{defn}
We recall the definition of subdifferentials. 
\begin{defn}
(Subdifferentials) Consider a function $f:\mathbb{R}^{n}\to\mathbb{R}\cup\{\infty\}$
at a point $\bar{x}$ where $f(\bar{x})$ is finite. Then the \emph{limiting
}(or \emph{Mordukhovich}) \emph{subdifferential} $\partial f(\bar{x})$,
\emph{horizon subdifferential} $\partial^{\infty}f(\bar{x})$ and
the \emph{Clarke (or generalized) subdifferential} $\partial_{C}f(\bar{x})$
are defined respectively by \begin{eqnarray*}
\partial f(\bar{x}) & := & \big\{v\mid(v,-1)\in N_{\scriptsize\epi(f)}\big(\bar{x},f(\bar{x})\big)\big\}\\
 & = & D^{*}E_{f}(\bar{x}\mid f(\bar{x}))(1),\\
\partial^{\infty}f(\bar{x}) & := & \big\{v\mid(v,0)\in N_{\scriptsize\epi(f)}\big(\bar{x},f(\bar{x})\big)\big\}\\
 & = & D^{*}E_{f}(\bar{x}\mid f(\bar{x}))(0),\\
\mbox{and }\partial_{C}f(\bar{x}) & := & \overline{\co}\partial f(\bar{x})\\
 & = & \overline{\co}D^{*}E_{f}(\bar{x}\mid f(\bar{x}))(1).\end{eqnarray*}

\end{defn}
The limiting and Clarke subdifferentials coincide with the usual definition
of subdifferential when $f$ is convex. The subdifferential $\partial f(\bar{x})$
gives important information on how $f$ varies with respect to $x$
when close to $\bar{x}$. 

We now recall the definition of generalized derivatives of set-valued
maps in the sense of \cite{set_diff}. Let $\mathbb{B}$ denote the
unit ball in the appropriate space.
\begin{defn}
\label{def:T-diff}\cite{set_diff} (Generalized differentiability)
Let $S:\mathbb{R}^{n}\rightrightarrows\mathbb{R}^{m}$ be such that
$S$ is locally closed at $(\bar{x},\bar{y})\in\gph(S)$, and let
$H:\mathbb{R}^{n}\rightrightarrows\mathbb{R}^{m}$ be a positively
homogeneous map. The map $S$ is \emph{pseudo strictly $H$-differentiable
at $(\bar{x},\bar{y})$ }if for any $\delta>0$, there are neighborhoods
$U_{\delta}$ of $\bar{x}$ and $V_{\delta}$ of $\bar{y}$ such that
\[
S(x)\cap V_{\delta}\subset S(x^{\prime})+H(x-x^{\prime})+\delta\|x-x^{\prime}\|\mathbb{B}\mbox{ for all }x,x^{\prime}\in U_{\delta}.\]
We shall also write \[
(H+\delta)(w):=H(w)+\delta\|w\|\mathbb{B}\]
to reduce notation. The map $S$ has the \emph{Aubin property} (or
the \emph{pseudo-Lipschitz }property) with modulus $\kappa\geq0$
if $S$ is pseudo strictly $H$-differentiable for some $H$ defined
by $H(w)=\kappa\|w\|\mathbb{B}$. The \emph{graphical modulus} is
the infimum of all such $\kappa$, and is denoted by $\lip\, S(\bar{x}\mid\bar{y})$.
\end{defn}
We now recall the definition of prefans and the generalized derivative
set $\mathcal{H}(D)$.
\begin{defn}
\cite{Iof81} (Prefans) We say that $H:\mathbb{R}^{n}\rightrightarrows\mathbb{R}^{m}$
is a \emph{prefan} if 
\begin{enumerate}
\item $H(p)$ is nonempty, convex and compact for all $p\in\mathbb{R}^{n}$,
\item $H$ is positively homogeneous, and
\item $\|H\|^{+}:=\sup_{\|w\|\leq1}\sup_{z\in H(w)}\|z\|$ is finite.\end{enumerate}
\begin{defn}
\cite{chardiff} (Generalized derivative set) Let $D:\mathbb{R}^{m}\rightrightarrows\mathbb{R}^{n}$
be a positively homogeneous, osc set-valued map s.t. $\|D\|^{+}$
is finite. We define the \emph{generalized derivative set }by \begin{eqnarray*}
\mathcal{H}(D) & := & \{H:\mathbb{R}^{n}\rightrightarrows\mathbb{R}^{m}:H\mbox{ is a prefan,}\\
 &  & \qquad\mbox{and for all }p\in\mathbb{R}^{n}\backslash\{0\}\mbox{ and }u\in\mathbb{R}^{m},\\
 &  & \qquad\underset{y\in H(p)}{\min}\left\langle u,y\right\rangle \leq\underset{v\in\scriptsize{\overline{\co}\,}D(u)}{\min}\left\langle v,p\right\rangle \}.\end{eqnarray*}

\end{defn}
\end{defn}
The Aubin criterion characterizes the graphical modulus $\lip S(\bar{x}\mid\bar{y})$
in terms of graphical derivatives (which are in turn defined in terms
of tangent cones), while the Mordukhovich criterion characterizes
$\lip S(\bar{x}\mid\bar{y})$ in terms of coderivatives.  Theorem
\ref{thm:boris-crit} and Lemma \ref{lem:conv-coderv-gen-derv} below
characterize the set of possible generalized derivatives at a point
$(\bar{x},\bar{y})\in\gph(S)$, and can be seen as a generalization
of the Mordukhovich criterion. While the proof in \cite{chardiff}
makes heavy use of graphical derivatives and recent work in \cite{DQZ06}
(who in turn acknowledged Frankowska's contribution), the main results
in finite dimensions have an appealing formulation in terms of coderivatives.
\begin{thm}
\cite{chardiff} \label{thm:boris-crit}(Characterization of generalized
derivatives) Let $S:\mathbb{R}^{n}\rightrightarrows\mathbb{R}^{m}$
be locally closed at $(\bar{x},\bar{y})\in\gph(S)$ and let $H:\mathbb{R}^{n}\rightrightarrows\mathbb{R}^{m}$
be a prefan. Then $S$ is pseudo strictly $H$-differentiable at $(\bar{x},\bar{y})$
if and only if $H\in\mathcal{H}\big(D^{*}S(\bar{x}\mid\bar{y})\big)$.
(Note that $\mathcal{H}\big(D^{*}S(\bar{x}\mid\bar{y})\big)=\mathcal{H}\big(\overline{\co}D^{*}S(\bar{x}\mid\bar{y})\big)$.)\end{thm}
\begin{lem}
\cite{chardiff}\label{lem:conv-coderv-gen-derv} (Convexified coderivatives
and generalized derivatives) Suppose $D_{i}:\mathbb{R}^{m}\rightrightarrows\mathbb{R}^{n}$
are positively homogeneous, osc, and $\|D_{i}\|^{+}$ are finite for
$i=1,2$. Then the following strict reverse inclusion properties hold:
\begin{enumerate}
\item $\mathcal{H}(D_{1})\supset\mathcal{H}(D_{2})$ iff $\overline{\co}\, D_{1}\subset\overline{\co}\, D_{2}$.
\item $\mathcal{H}(D_{1})\supsetneq\mathcal{H}(D_{2})$ iff $\overline{\co}\, D_{1}\subsetneq\overline{\co}\, D_{2}$.
\item $\mathcal{H}(D_{1})=\mathcal{H}(D_{2})$ iff $\overline{\co}\, D_{1}=\overline{\co}\, D_{2}$.
\end{enumerate}
\end{lem}
These results show that the convexified limiting coderivative $\overline{\co}D^{*}S(\cdot|\cdot)(\cdot)$
is an effective tool for studying the generalized derivatives of set-valued
maps, just like the way the Clarke subdifferential is useful for studying
the generalized differentiability of single-valued maps.

We recall the definition of inner semicompactness that will be used
in the chain rules for set-valued maps in this paper.
\begin{defn}
(Inner semicompactness) We say that $S:\mathbb{R}^{n}\rightrightarrows\mathbb{R}^{m}$
is \emph{inner semicompact} at $\bar{x}\in\dom(S)$ if for every sequence
$x_{k}\to\bar{x}$, there is a sequence $y_{k}\in S(x_{k})$ that
contains a convergent subsequence as $k\to\infty$. 
\end{defn}
In finite dimensions, if there is a neighborhood $U$ of $\bar{x}$
and a bounded neighborhood $V$ such that $S(U)\subset V$, then $S$
is inner semicompact at $\bar{x}$.

Finally, we recall the definition of regularity and a straightforward
consequence of graphical regularity.
\begin{defn}
(Regularity) We say that $C\subset\mathbb{R}^{n}$ is \emph{Clarke
regular }at $\bar{x}\in C$ if $C$ is locally closed at $\bar{x}$
and $N_{C}(\bar{x})=\hat{N}_{C}(\bar{x})$. We say that $S:\mathbb{R}^{n}\rightrightarrows\mathbb{R}^{m}$
is \emph{graphically regular} at $(\bar{x},\bar{y})\in\gph(S)$ if
$\gph(S)$ is Clarke regular at $(\bar{x},\bar{y})$.\end{defn}
\begin{fact}
(Convexified limiting coderivatives under graph regularity) If $S:\mathbb{R}^{n}\rightrightarrows\mathbb{R}^{m}$
is graphically regular\emph{ }at $(\bar{x},\bar{y})\in\gph(S)$, then
$\gph(D^{*}S(\bar{x}\mid\bar{y}))=\gph(\hat{D}^{*}S(\bar{x}\mid\bar{y}))$.
Furthermore, $\gph(\hat{D}^{*}S(\bar{x}\mid\bar{y}))$ is a convex
cone, and we have $\overline{\co}D^{*}S(\bar{x}\mid\bar{y})\equiv D^{*}S(\bar{x}\mid\bar{y})$.
\end{fact}

\section{\label{sec:calculus-co-coderv}Calculus of convexified limiting coderivatives}

In this section, we discuss how the chain rule for the convexified
limiting coderivatives can be obtained directly from the  coderivative
chain rules, removing parts irrelevant in the finite dimensional case.
In Lemma \ref{lem:gen-opt}, we deduce that the convexified limiting
coderivative, together with the limiting subdifferential, are sufficient
in calculating the Clarke subdifferential of marginal functions. This
suggests that the convexified limiting coderivative of the reachable
map as calculated in \ref{sec:cts-incl}, while not as precise as
the coderivative, can be a satisfactory conclusion. 

We first write down the chain rule for finite dimensional coderivatives
based on \cite[Theorem 3.13]{Mor06} and \cite[Theorem 10.37]{RW98}.
The formulas \eqref{eq:co-chain-rule-1} and \eqref{eq:co-chain-rule-2}
for convexified limiting coderivatives are straightforward. 
\begin{thm}
\label{thm:Coderivative-chain-rule}(Coderivative chain rule) Let
$G:\mathbb{R}^{l}\rightrightarrows\mathbb{R}^{m}$, $F:\mathbb{R}^{m}\rightrightarrows\mathbb{R}^{n}$,
$\bar{z}\in(F\circ G)(\bar{x})$, and \[
S(x,z):=G(x)\cap F^{-1}(z)=\{y\in G(x):z\in F(y)\}.\]
The following assertions hold:
\begin{enumerate}
\item Given $\bar{y}\in S(\bar{x},\bar{z})$, assume that $S$ is inner
semicontinuous at $(\bar{x},\bar{z},\bar{y})$, that the graphs of
$F$ and $G$ are locally closed around the points $(\bar{y},\bar{z})$
and $(\bar{x},\bar{y})$ respectively, and that the qualification
condition \begin{equation}
D^{*}F(\bar{y}\mid\bar{z})(0)\cap-D^{*}G^{-1}(\bar{y}\mid\bar{x})(0)=\{0\}\label{eq:chain-rule-CQ}\end{equation}
is fulfilled. Then one has \[
D^{*}(F\circ G)(\bar{x}\mid\bar{z})\subset D^{*}G(\bar{x}\mid\bar{y})\circ D^{*}F(\bar{y}\mid\bar{z}),\]
which in turn implies \begin{equation}
\overline{\co}D^{*}(F\circ G)(\bar{x}\mid\bar{z})\subset\overline{\co}D^{*}G(\bar{x}\mid\bar{y})\circ D^{*}F(\bar{y}\mid\bar{z}).\label{eq:co-chain-rule-1}\end{equation}
 
\item Assume that $S$ is inner semicompact at $(\bar{x},\bar{z})$, that
$G$ and $F^{-1}$ are closed-graph whenever $x$ is near $\bar{x}$
and $z$ is near $\bar{z}$, respectively, and that \eqref{eq:chain-rule-CQ}
holds for every $\bar{y}\in S(\bar{x},\bar{z})$. Then \[
D^{*}(F\circ G)(\bar{x}\mid\bar{z})\subset\bigcup_{\bar{y}\in S(\bar{x},\bar{z})}D^{*}G(\bar{x}\mid\bar{y})\circ D^{*}F(\bar{y}\mid\bar{z}),\]
which in turn implies \begin{equation}
\overline{\co}D^{*}(F\circ G)(\bar{x}\mid\bar{z})\subset\overline{\co}\bigcup_{\bar{y}\in S(\bar{x},\bar{z})}\overline{\co}D^{*}G(\bar{x}\mid\bar{y})\circ D^{*}F(\bar{y}\mid\bar{z}).\label{eq:co-chain-rule-2}\end{equation}

\item If $S$ is locally bounded at $(\bar{x},\bar{z})$, \eqref{eq:chain-rule-CQ}
holds for every $\bar{y}\in S(\bar{x},\bar{z})$, and $F$ and $G$
are both graph convex (i.e., have convex graphs), then $F\circ G$
is also graph convex, and \[
D^{*}(F\circ G)(\bar{x}\mid\bar{z})=D^{*}G(\bar{x}\mid\bar{y})\circ D^{*}F(\bar{y}\mid\bar{z})\mbox{ for any }\bar{y}\in S(\bar{x},\bar{z}).\]

\end{enumerate}
\end{thm}
The formula \eqref{eq:co-chain-rule-2} is not any stronger if its
RHS is replaced by \[
\overline{\co}\bigcup_{\bar{y}\in S(\bar{x},\bar{z})}D^{*}G(\bar{x}\mid\bar{y})\circ D^{*}F(\bar{y}\mid\bar{z}),\]
since this formula is equal to the RHS of \eqref{eq:co-chain-rule-2}.
Therefore, to find the convexified limiting coderivative $\overline{\co}D^{*}(F\circ G)(\bar{x}\mid\bar{z})$,
the convexified limiting coderivative of $G$, i.e., $\overline{\co}D^{*}G$,
is sufficient. We explore the possibilities if we had relaxed the
formulas \eqref{eq:co-chain-rule-1} and \eqref{eq:co-chain-rule-2}
by replacing the relevant formulas with $\overline{\co}D^{*}G(\bar{x}\mid\bar{y})\circ\overline{\co}D^{*}F(\bar{y}\mid\bar{z})$
instead.
\begin{example}
\label{exa:chain-rule-eg}(Tightness of chain rules) Consider the
set-valued maps $G_{i}:\mathbb{R}\rightrightarrows\mathbb{R}$, $i=1,2,3$
and $F:\mathbb{R}\rightrightarrows\mathbb{R}$ defined by \begin{eqnarray*}
G_{1}(x) & = & \begin{cases}
\mathbb{R} & \mbox{ if }x\leq0\\
(-\infty,-x]\cup[x,\infty) & \mbox{ if }x\geq0,\end{cases}\\
G_{2}(x) & = & \big[\min(0,x),\infty\big),\\
G_{3}(x) & = & \big[\max(x/2,x),\infty\big),\\
\mbox{ and }F(x) & = & \{-x\}\cup\{x\}.\end{eqnarray*}
As illustrated in Table \ref{tab:chain-rule-eg}, we have\begin{eqnarray*}
\overline{\co}D^{*}(F\circ G_{1})(0\mid0)=\overline{\co}D^{*}G_{1}(0\mid0)\circ D^{*}F(0\mid0)\subsetneq\overline{\co}D^{*}G_{1}(0\mid0)\circ\overline{\co}D^{*}F(0\mid0),\\
\overline{\co}D^{*}(F\circ G_{2})(0\mid0)\subsetneq\overline{\co}D^{*}G_{2}(0\mid0)\circ D^{*}F(0\mid0)=\overline{\co}D^{*}G_{2}(0\mid0)\circ\overline{\co}D^{*}F(0\mid0),\\
\mbox{and }\overline{\co}D^{*}(F\circ G_{3})(0\mid0)\subsetneq\overline{\co}D^{*}G_{3}(0\mid0)\circ D^{*}F(0\mid0)\subsetneq\overline{\co}D^{*}G_{3}(0\mid0)\circ\overline{\co}D^{*}F(0\mid0).\end{eqnarray*}

\end{example}
\begin{table}
\begin{tabular}{|r|c|c|c|}
\hline 
 & $i=1$ & $2$ & $3$\tabularnewline
\hline
\hline 
$G_{i}$ & \includegraphics[scale=0.25]{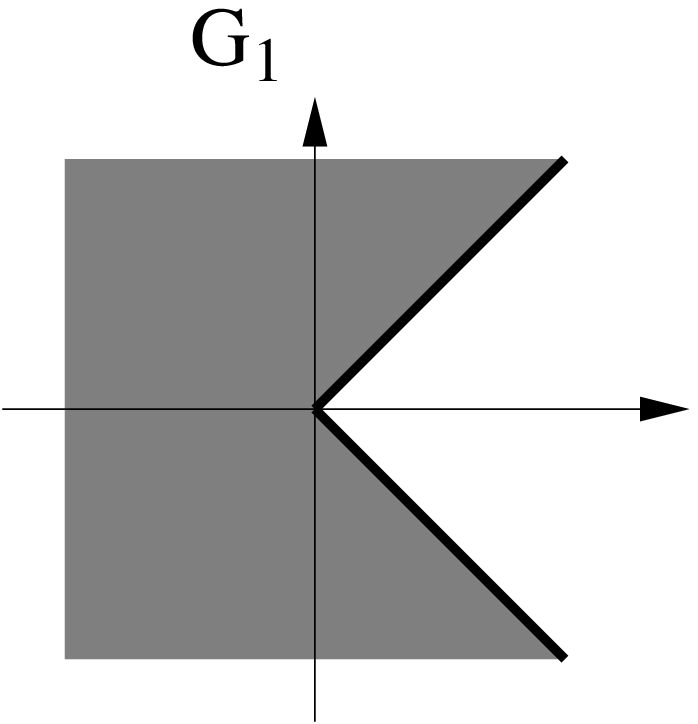} & \includegraphics[scale=0.25]{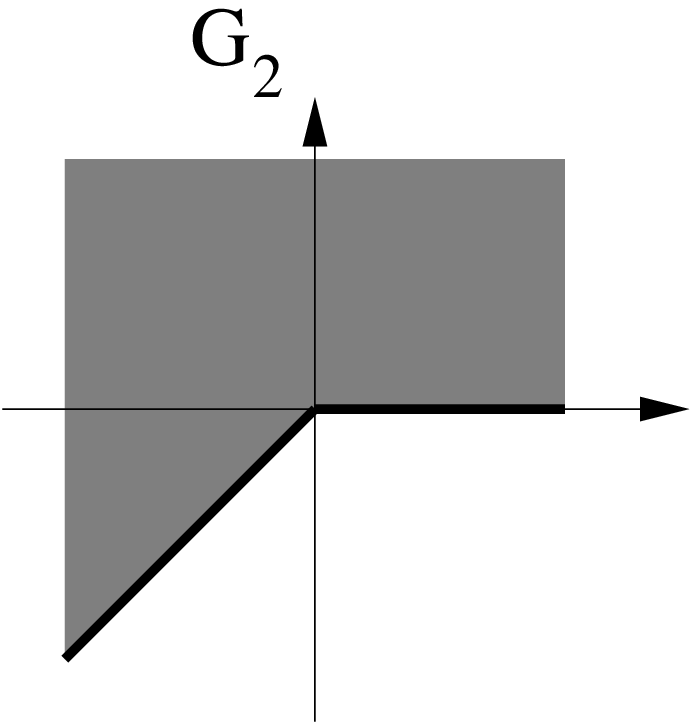} & \includegraphics[scale=0.25]{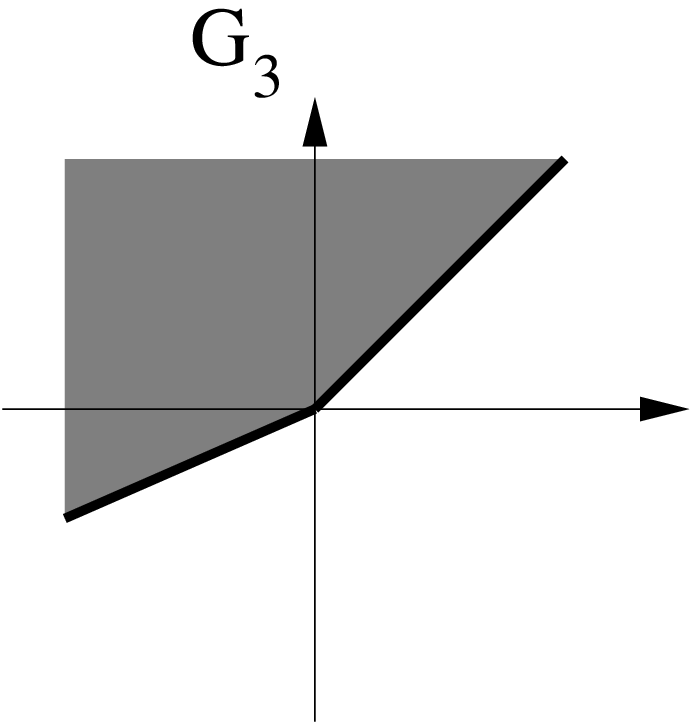}\tabularnewline
\hline 
$F$ & \includegraphics[scale=0.25]{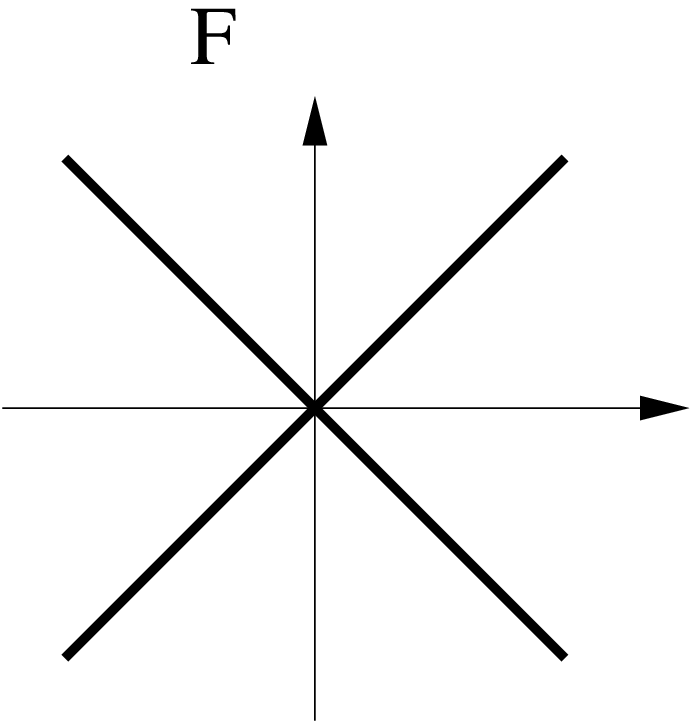} & \includegraphics[scale=0.25]{Ffn} & \includegraphics[scale=0.25]{Ffn}\tabularnewline
\hline 
$F\circ G_{i}$ & \includegraphics[scale=0.25]{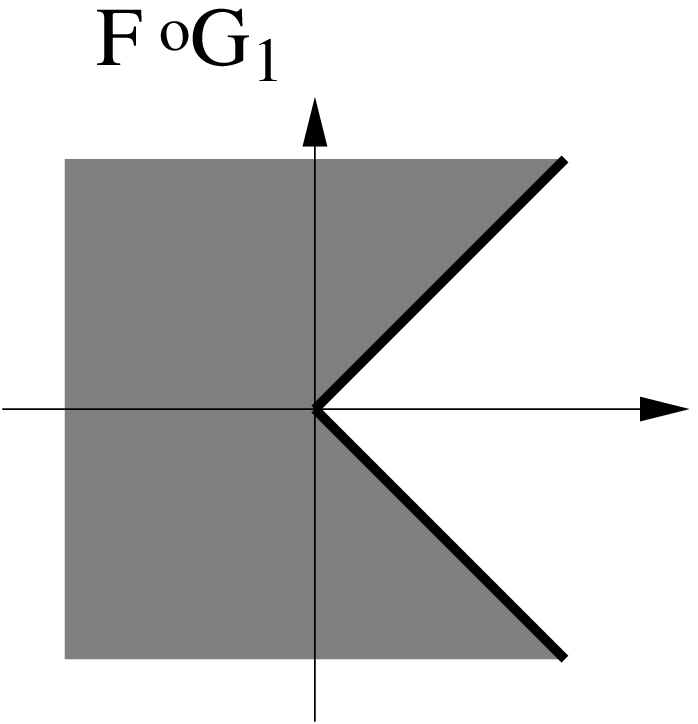} & \includegraphics[scale=0.25]{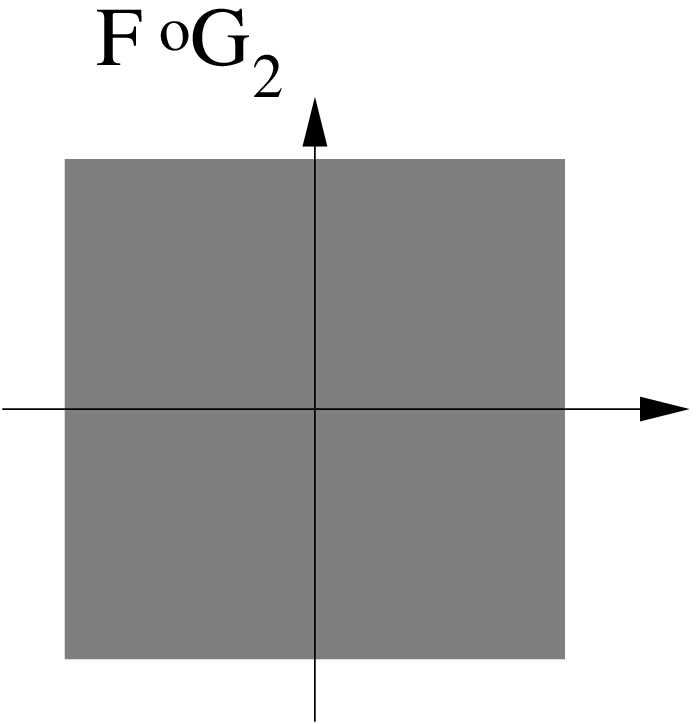} & \includegraphics[scale=0.25]{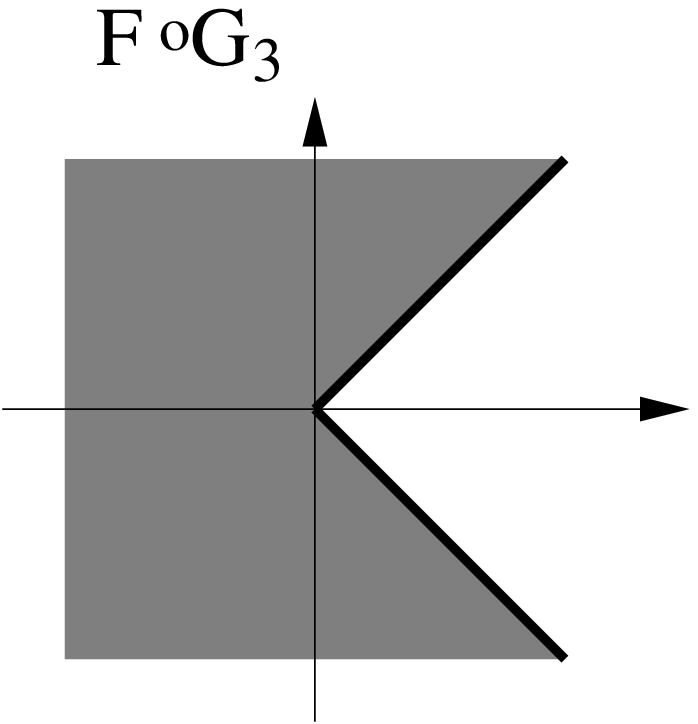}\tabularnewline
\hline 
$D^{*}F(0\mid0)$ & \includegraphics[scale=0.25]{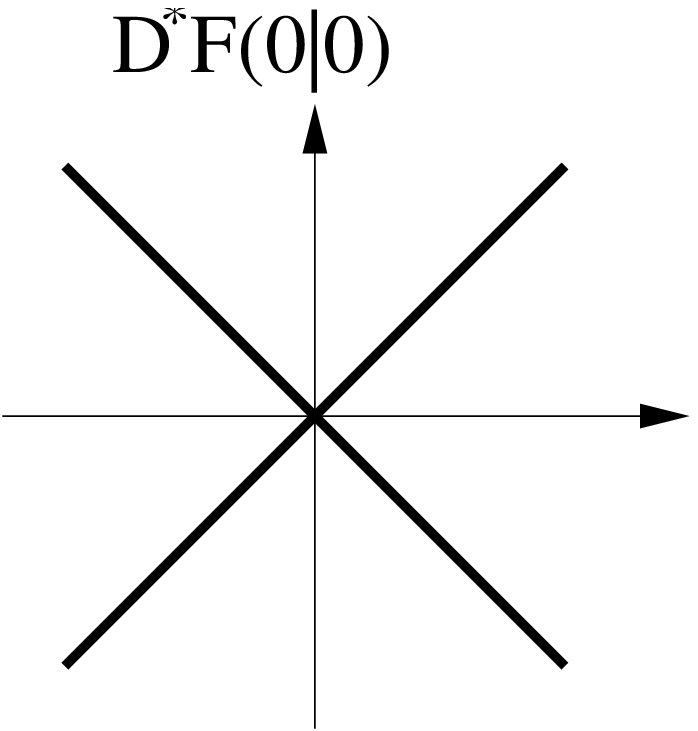} & \includegraphics[scale=0.25]{DF} & \includegraphics[scale=0.25]{DF}\tabularnewline
\hline 
$D^{*}G_{i}(0\mid0)$ & \includegraphics[scale=0.25]{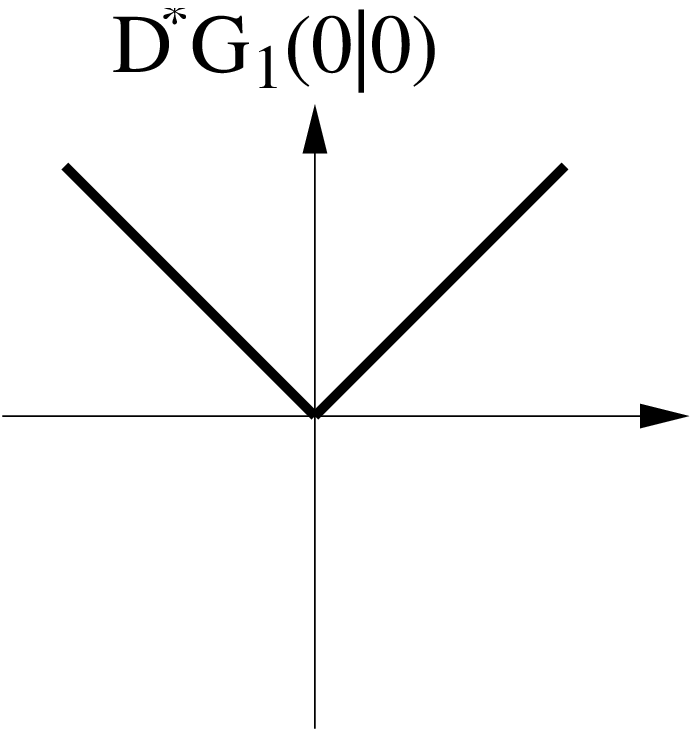} & \includegraphics[scale=0.25]{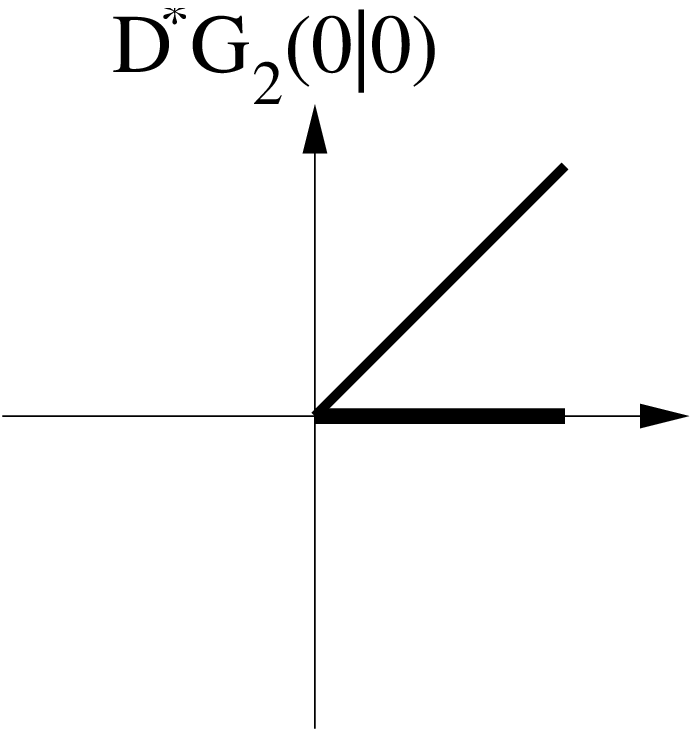} & \includegraphics[scale=0.25]{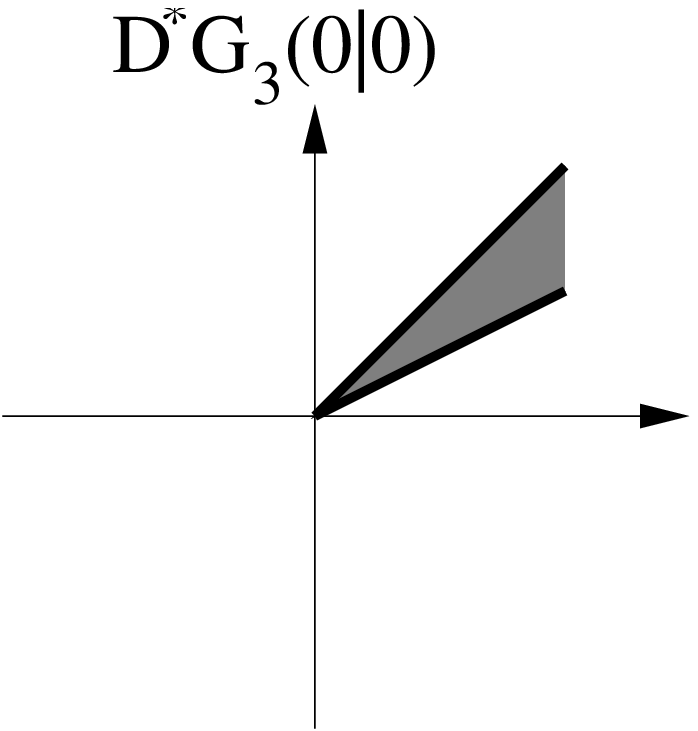}\tabularnewline
\hline 
$\overline{\co}D^{*}(F\circ G_{i})(0\mid0)$ & \includegraphics[scale=0.25]{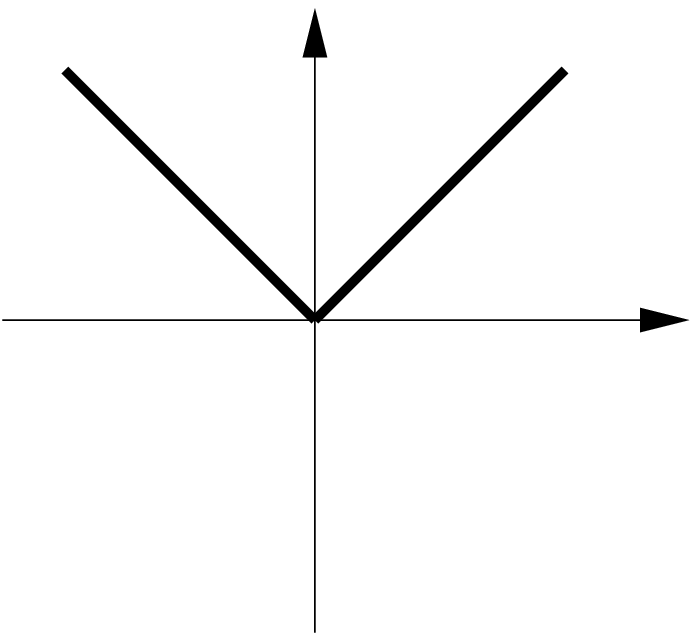} & \includegraphics[scale=0.25]{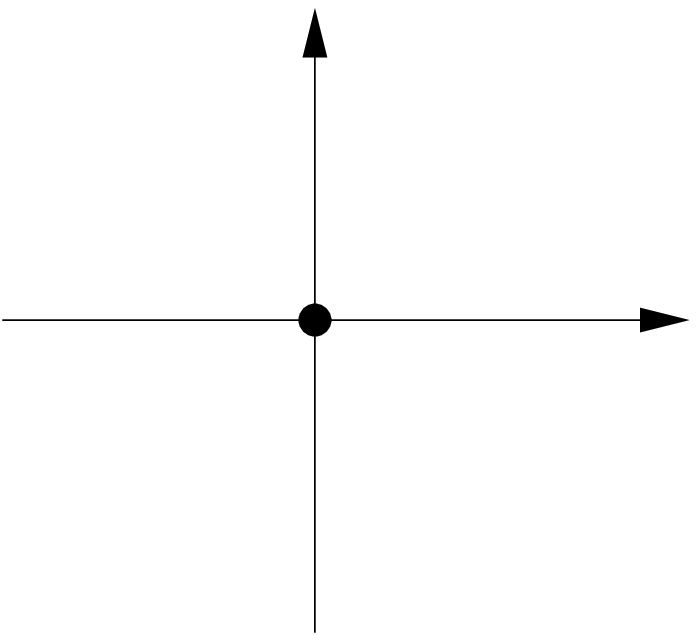} & \includegraphics[scale=0.25]{DFG1}\tabularnewline
\hline 
$\overline{\co}D^{*}G_{i}(0\mid0)\circ D^{*}F(0\mid0)$ & \includegraphics[scale=0.25]{DFG1} & \includegraphics[scale=0.25]{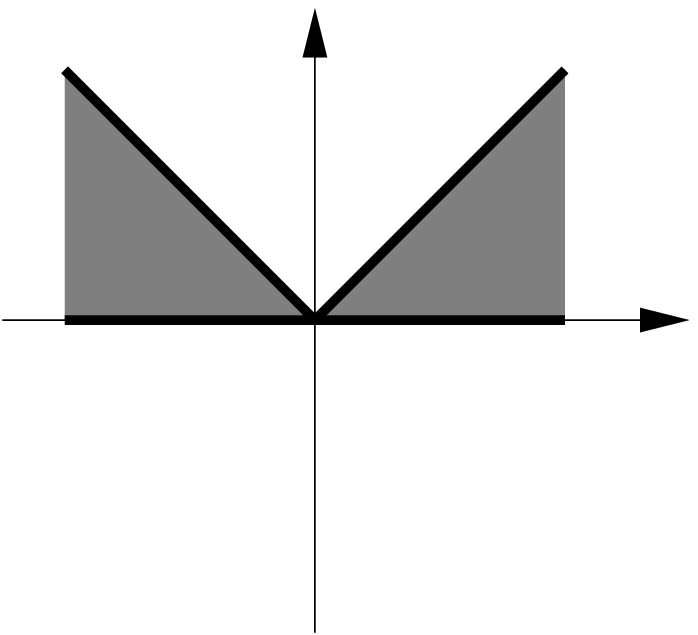} & \includegraphics[scale=0.25]{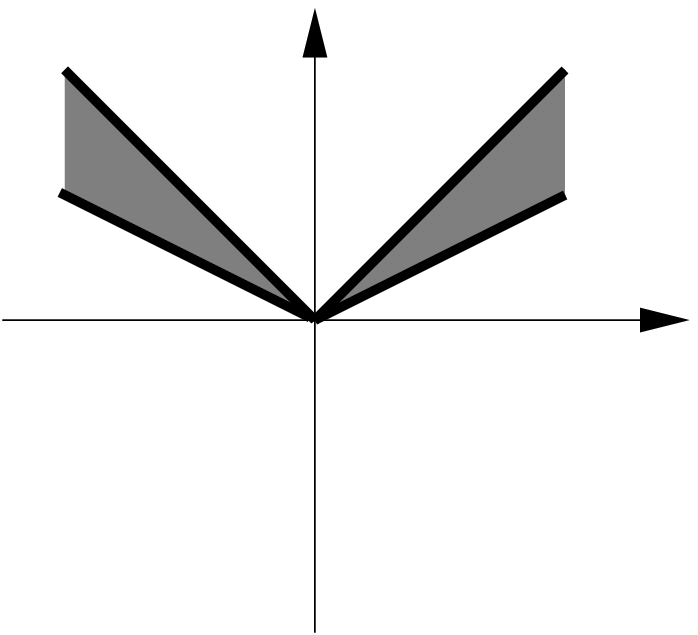}\tabularnewline
\hline 
$\overline{\co}D^{*}G_{i}(0\mid0)\circ\overline{\co}D^{*}F(0\mid0)$ & \includegraphics[scale=0.25]{DGcF3} & \includegraphics[scale=0.25]{DGcF3} & \includegraphics[scale=0.25]{DGcF3}\tabularnewline
\hline
\end{tabular}

\caption{\label{tab:chain-rule-eg}Possible scenarios in chain rule of set-valued
maps from Example \ref{exa:chain-rule-eg}.}

\end{table}

The following general principle in the optimization of marginal functions
will be used later. We take this result from \cite[Theorem 3.38]{Mor06}.
\begin{lem}
\label{lem:gen-opt}(Subdifferential of marginal functions) For the
marginal function \eqref{eq:marginal-1}, define the argminimum mapping
by \[
M(x):=\{y\in G(x)\mid\varphi(x,y)=f(x)\}.\]
The following hold: 
\begin{enumerate}
\item Given $\bar{y}\in M(\bar{x})$, assume that $M$ is inner semicontinuous
at $(\bar{x},\bar{y})$, that $\varphi(x,y)$ is l.s.c. around $(\bar{x},\bar{y})$,
and that $\gph(G)$ is locally closed at $(\bar{x},\bar{y})$. Suppose
also that the qualification condition \begin{equation}
\partial^{\infty}\varphi(\bar{x},\bar{y})\cap-N_{\scriptsize\gph(G)}(\bar{x},\bar{y})=\{0\}\label{eq:subdiff-chain-CQ}\end{equation}
is satisfied. Then one has the inclusion \begin{eqnarray}
\partial f(\bar{x}) & \subset & \bigcup_{(x^{*},y^{*})\in\partial\varphi(\bar{x},\bar{y})}[x^{*}+D^{*}G(\bar{x}\mid\bar{y})(y^{*})].\label{eq:subdiff-chain-incl}\\
\mbox{ and }\partial_{C}f(\bar{x}) & \subset & \overline{\co}\bigcup_{(x^{*},y^{*})\in\partial\varphi(\bar{x},\bar{y})}[x^{*}+\overline{\co}D^{*}G(\bar{x}\mid\bar{y})(y^{*})].\nonumber \end{eqnarray}

\item Assume that $M$ is inner semicompact at $\bar{x}$, that $G$ is
closed-graph and $\varphi$ is l.s.c. on $\gph(G)$ whenever $x$
is near $\bar{x}$, and that the other assumptions in (1) are satisfied
for every $\bar{y}\in M(\bar{x})$. Then one has analogs of inclusion
\eqref{eq:subdiff-chain-incl}, where the sets on the right-hand sides
are replaced by their unions over $\bar{y}\in M(\bar{x})$. 
\item Assume that $M$ is locally bounded at $\bar{x}$, \eqref{eq:subdiff-chain-CQ}
is satisfied for every $\bar{y}\in M(\bar{x})$, $G$ is graph-convex
and $\varphi$ is convex. Then $f$ is convex, and \begin{eqnarray*}
\partial f(\bar{x}) & = & \{x^{*}+D^{*}G(\bar{x}\mid\bar{y})(y^{*})\mid(x^{*},y^{*})\in\partial\varphi(\bar{x},\bar{y})\}\\
 &  & \mbox{for any }\bar{y}\in M(\bar{x}).\end{eqnarray*}

\end{enumerate}
\end{lem}
\begin{proof}
 Cases (1) and (2) are exactly the statement of \cite[Theorem 3.38]{Mor06},
and we prove only (3) from Theorem \ref{thm:Coderivative-chain-rule}(3).
Consider the map $\bar{G}:\mathbb{R}^{n}\rightrightarrows\mathbb{R}^{n}\times\mathbb{R}^{m}$
defined by $\bar{G}(x)=\{x\}\times G(x)$. The coderivative $D^{*}\bar{G}(\bar{x}\mid(\bar{x},\bar{y})):\mathbb{R}^{n}\times\mathbb{R}^{m}\rightrightarrows\mathbb{R}^{n}$
is easily evaluated to be \[
D^{*}\bar{G}\big(\bar{x}\mid(\bar{x},\bar{y})\big)(p,q)=p+D^{*}G(\bar{x}\mid\bar{y})(q).\]
Noting that $E_{f}=E_{\varphi}\circ\bar{G}$, the constraint qualification
we need to check in Theorem \ref{thm:Coderivative-chain-rule}(3)
is \begin{equation}
D^{*}E_{\varphi}\big((\bar{x},\bar{y})\mid f(\bar{x})\big)(0)\cap-D^{*}\bar{G}^{-1}\big((\bar{x},\bar{y})\mid\bar{x}\big)(0)=\{0\}.\label{eq:convex-subdif-chain-to-check}\end{equation}
Note that $\partial^{\infty}\varphi(\bar{x},\bar{y})=D^{*}E_{\varphi}((\bar{x},\bar{y})\mid f(\bar{x}))(0)$.
Now, $(p,q)\in D^{*}\bar{G}^{-1}((\bar{x},\bar{y})\mid\bar{x})(0)$
if and only if $(p,q,0)\in N_{\scriptsize\gph(\bar{G}^{-1})}(\bar{x},\bar{y},\bar{x})$,
which is in turn equivalent to $(0,p,q)\in N_{\scriptsize\gph(\bar{G})}(\bar{x},\bar{x},\bar{y})$.
We see that $\gph(\bar{G})$ is the image of a linear map of $\gph(G)$
and use a rule of normal cones on linear maps in \cite[Theorem 6.43]{RW98}
to obtain \[
N_{\scriptsize\gph(\bar{G})}(\bar{x},\bar{x},\bar{y})=\{(u,w,v)\mid(u+w,v)\in N_{\scriptsize\gph(G)}(\bar{x},\bar{y})\}.\]
Thus $(0,p,q)\in N_{\scriptsize\gph(\bar{G})}(\bar{x},\bar{x},\bar{y})$
iff $(p,q)\in N_{\scriptsize\gph(G)}(\bar{x},\bar{y})$. Therefore
\eqref{eq:convex-subdif-chain-to-check} is equivalent to \eqref{eq:subdiff-chain-CQ}.
We then apply Theorem \ref{thm:Coderivative-chain-rule}(3) to get
\begin{eqnarray*}
\partial f(\bar{x}) & = & D^{*}E_{f}\big(\bar{x}\mid f(\bar{x})\big)(1)\\
 & = & D^{*}(E_{\varphi}\circ\bar{G})\big(\bar{x}\mid f(\bar{x})\big)(1).\end{eqnarray*}
Then for any $\bar{y}\in M(\bar{x})$, \begin{eqnarray*}
\partial f(\bar{x}) & = & D^{*}\bar{G}\big(\bar{x}\mid(\bar{x},\bar{y})\big)\circ D^{*}E_{\varphi}\big((\bar{x},\bar{y})\mid f(\bar{x})\big)(1)\\
 & = & D^{*}\bar{G}\big(\bar{x}\mid(\bar{x},\bar{y})\big)\big(\partial\varphi(\bar{x},\bar{y})\big)\\
 & = & \{x^{*}+D^{*}G(\bar{x}\mid\bar{y})(y^{*})\mid(x^{*},y^{*})\in\partial\varphi(\bar{x},\bar{y})\}.\end{eqnarray*}

\end{proof}
We remark that \cite[Section 10H]{RW98} and \cite[Section 3.2]{Mor06}
contain other coderivative calculus rules that can be easily extended
for the convexified limiting coderivative. As we have remarked after
Theorem \ref{thm:Coderivative-chain-rule}, the convexified limiting
coderivative of $G$ in Lemma \ref{lem:gen-opt} is sufficient for
obtaining the Clarke subdifferential of $f$.
\begin{rem}
\label{rem:alternate-marginal}(Alternative view of marginal functions)
A different view useful for later discussions is to consider \begin{eqnarray*}
 & \underset{\scriptsize{(x,y)}}{\min} & \varphi(x,y)\\
 & \mbox{s.t.} & (x,y)\in\gph(G).\end{eqnarray*}
As is well known in nonlinear programming, if the point $(\bar{x},\bar{y})$
is optimal, then $0\in\partial\varphi(\bar{x},\bar{y})+N_{\scriptsize\gph(G)}(\bar{x},\bar{y})$.
Recall that through the definition of coderivatives, $N_{\scriptsize\gph(G)}(\bar{x},\bar{y})$
is related to $\gph(D^{*}G(\bar{x}\mid\bar{y}))$ by a linear transformation.
\end{rem}

\section{\label{sec:discrete-incl}Subdifferential analysis of discretized
inclusions}

In this section, we consider the discretized inclusion and calculate
the coderivatives of its reachable map. One can then obtain the subdifferential
dependence of the differential inclusion in terms of its initial conditions.
We can then obtain a necessary optimality condition of the discretized
inclusion similar to the Euler-Lagrange and Transversality conditions.
Finally, we discuss the limitations of obtaining a discretized version
of the Weierstrass-Pontryagin maximum principle. 

We consider the following discrete inclusion as the analogue to the
differential inclusion \eqref{eq:main-diff-in}:\begin{eqnarray}
 & \underset{x_{k}\in\mathbb{R}^{n}\scriptsize{\mbox{ for }}k\in\overline{0,N}}{\min} & \varphi(x_{0},x_{N})\label{eq:main-discrete-in}\\
 & \mbox{s.t.} & x_{k}\in x_{k-1}+(\Delta t)F\big((k-1)(\Delta t),x_{k-1}\big).\nonumber \end{eqnarray}
Here, $\Delta t=T/N$. The inclusion systems above can be further
modified to one defined in terms of the reachable map. The discretized
version of the reachable map $R_{N}:\mathbb{R}^{n}\rightrightarrows\mathbb{R}^{n}$
can be defined by \begin{eqnarray}
R_{N}(x_{0}) & = & \{x_{N}:\exists x_{k}\in\mathbb{R}^{n}\mbox{ for }k=\overline{1,N-1}\mbox{ s.t. }\label{eq:R-N-discrete-reachable}\\
 &  & \qquad x_{k}\in x_{k-1}+(\Delta t)F\big((k-1)(\Delta t),x_{k-1}\big)\mbox{ for all }k=\overline{1,N}\}.\nonumber \end{eqnarray}
Then \eqref{eq:main-discrete-in} can be rewritten as \begin{eqnarray}
 & \underset{x_{0},x_{N}}{\min} & \varphi(x_{0},x_{N})\label{eq:reachable-discrete-in}\\
 & \mbox{s.t.} & x_{N}\in R_{N}(x_{0})\subset\mathbb{R}^{n}.\nonumber \end{eqnarray}

\begin{thm}
\label{thm:discrete-coderv}(Coderivatives of discretized reachable
map) Recall the map $R_{N}:\mathbb{R}^{n}\rightrightarrows\mathbb{R}^{n}$
as defined in \eqref{eq:reachable-discrete-in}. Let $\Delta t=T/N$,
and define $F_{k,N}:\mathbb{R}^{n}\rightrightarrows\mathbb{R}^{n}$
and $M_{k,N}:\mathbb{R}^{n}\rightrightarrows\mathbb{R}^{n}$ by \begin{equation}
F_{k,N}(\cdot):=F(k(\Delta t),\cdot)\qquad\mbox{and }\qquad M_{k,N}(x):=x+(\Delta t)F_{k,N}(x).\label{eq:FkN-and-MkN}\end{equation}
Note that $R_{N}=M_{N-1,N}\circ M_{N-2,N}\circ\cdots\circ M_{0,N}$.
Assume that $F_{k,N}$ is locally Lipschitz and for all $x$, $k$
and $N$, $F_{k,N}(\cdot)$ is locally bounded at $x$, i.e., there
exists a neighborhood $U$ of $x$ and finite $R$ such that $F_{k,N}(x^{\prime})\subset R\mathbb{B}$
for all $x^{\prime}\in U$.
\begin{enumerate}
\item For $x_{N}\in R(x_{0})$, the coderivative of $R_{N}:\mathbb{R}^{n}\rightrightarrows\mathbb{R}^{n}$
satisfies\begin{eqnarray}
 &  & D^{*}R_{N}(x_{0}\mid x_{N})\nonumber \\
 & \subset & \bigcup_{\{\tilde{x}_{i}\}_{i=0}^{N}\in\mathcal{X}_{N}}\underbrace{D^{*}M_{0,N}(\tilde{x}_{0}\mid\tilde{x}_{1})\circ\cdots\circ D^{*}M_{N-1,N}(\tilde{x}_{N-1}\mid\tilde{x}_{N})}_{G_{\{\tilde{x_{i}}\}_{i=0}^{N}}}\label{eq:D-R-N-first-line}\end{eqnarray}
where \begin{equation}
\mathcal{X}_{N}=\big\{\{\tilde{x}_{i}\}_{i=0}^{N}:\tilde{x}_{k}\in M_{k-1,N}(\tilde{x}_{k-1})\mbox{ for all }k\in\overline{1,N}\mbox{, }\tilde{x}_{0}=x_{0}\mbox{ and }\tilde{x}_{N}=x_{N}\big\}.\label{eq:X-N}\end{equation}
\textup{\emph{ }}
\item If in addition $F_{k,N}(\cdot)$ are all graph convex, then \[
D^{*}R_{N}(x_{0}\mid x_{N})=D^{*}M_{0,N}(\tilde{x}_{0}\mid\tilde{x}_{1})\circ\cdots\circ D^{*}M_{N-1,N}(\tilde{x}_{N-1}\mid\tilde{x}_{N})\]
 for any $\{\tilde{x}_{i}\}_{i=0}^{N}\in\mathcal{X}_{N}$. 
\item \textup{\emph{Consider the conditions:}}

\begin{enumerate}
\item \textup{\emph{$p_{0}\in D^{*}R_{N}(x_{0}|x_{N})(p_{N})$ }}
\item \textup{\emph{There are $\{\tilde{x}_{i}\}_{i=0}^{N}\in\mathcal{X}_{N}$
and $\{\tilde{p}_{i}\}_{i=0}^{N}$ such that $p_{0}=\tilde{p}_{0}$,
$p_{N}=\tilde{p}_{N}$ and \begin{equation}
\frac{p_{k}-p_{k-1}}{\Delta t}\in-D^{*}F_{k-1,N}\left(\tilde{x}_{k-1}\mid\frac{1}{\Delta t}(\tilde{x}_{k}-\tilde{x}_{k-1})\right)(p_{k})\mbox{ for all }k\in\overline{1,N}.\label{eq:p-k-second-line}\end{equation}
}}
\end{enumerate}
\end{enumerate}
\textup{\emph{We have (a) implies (b), and in the case where each
$F_{k,N}(\cdot)$ is graph convex for all $k\in\overline{0,(N-1)}$,
the converse holds as well.}}

\end{thm}
\begin{proof}
 For (1), the case where $N=2$ follows directly from Theorem \ref{thm:Coderivative-chain-rule}(2).
The local Lipschitz continuity of $F_{k,N}(\cdot)$ implies that the
graph of $F_{k,N}(\cdot)$ is closed. The local boundedness of $F_{k,N}$
ensures that $S(\cdot,\cdot)$ in Theorem \ref{thm:Coderivative-chain-rule}(2)
is inner semicompact, and the local Lipschitz continuity implies the
constraint qualification in \eqref{eq:chain-rule-CQ} holds. The
case for general $N$ is easily deduced from the case where $N=2$.
For (2), we follow the similar steps and apply Theorem \ref{thm:Coderivative-chain-rule}(3).

 To prove that (3a) implies (3b), \emph{let $G_{\{\tilde{x}_{i}\}_{i=0}^{N}}:\mathbb{R}^{n}\rightrightarrows\mathbb{R}^{n}$
be the formula as marked in \eqref{eq:D-R-N-first-line}. For $p_{N},p_{0}\in\mathbb{R}^{n}$,
we have $p_{0}\in G_{\{\tilde{x}_{i}\}_{i=0}^{N}}(p_{N})$ if and
only if there exists some $\{\tilde{p}_{i}\}_{i=0}^{N}$ such that
$p_{0}=\tilde{p}_{0}$, $p_{N}=\tilde{p}_{N}$ and \begin{equation}
p_{k-1}\in D^{*}M_{k-1}(\tilde{x}_{k-1}\mid\tilde{x}_{k})(p_{k})\mbox{ for all }k\in\overline{1,N},\label{eq:p-k-first-line}\end{equation}
}From the definition of $M_{k-1,N}$ and calculus rules for coderivatives
in \cite[Section 10H]{RW98}, we have \begin{equation}
D^{*}M_{k-1,N}(\tilde{x}_{k-1}\mid\tilde{x}_{k})=I+(\Delta t)D^{*}F_{k-1,N}\left(\tilde{x}_{k-1}\mid\frac{1}{\Delta t}(\tilde{x}_{k}-\tilde{x}_{k-1})\right).\label{eq:coderiv-M-k-N}\end{equation}
The formula \eqref{eq:p-k-second-line} follows easily from (1).
The converse holds due to (2).
\end{proof}
Putting together the previous results, we have the following necessary
optimality condition for the discrete inclusion problem.
\begin{thm}
\label{thm:nec-condn-discrete}(Subdifferential analysis of discrete
inclusions) For the discrete inclusion \eqref{eq:main-discrete-in},
suppose $F(t,\cdot)$ is Lipschitz and for each $t$, there is some
$b(t)<\infty$ such that $F(t,x)\subset b(t)\mathbb{B}$ for all $x$.
Define $F_{k,N}:\mathbb{R}^{n}\rightrightarrows\mathbb{R}^{n}$ and
$M_{k,N}:\mathbb{R}^{n}\rightrightarrows\mathbb{R}^{n}$ as in \eqref{eq:FkN-and-MkN},
$R_{N}:\mathbb{R}^{n}\rightrightarrows\mathbb{R}^{n}$ by \eqref{eq:R-N-discrete-reachable},
and $f:\mathbb{R}^{n}\to\mathbb{R}$ by \begin{eqnarray*}
f(x_{0}):= & \underset{x_{N}}{\min} & \varphi(x_{0},x_{N})\\
 & \mbox{s.t.} & x_{N}\in R_{N}(x_{0})\subset\mathbb{R}^{n}.\end{eqnarray*}
Suppose $\varphi(\cdot,\cdot)$ is lsc. Then \[
\partial f(x_{0})\subset\bigcup_{\scriptsize{\begin{array}{c}
\{\tilde{x}_{i}\}_{i=0}^{N}\in\mathcal{X}_{N},x_{N}\in R_{N}(x_{0})\\
\varphi(x_{0},x_{N})=f(x_{0})\end{array}}}\{x^{*}+D^{*}G_{\{\tilde{x}_{i}\}_{i=0}^{N}}(x_{0}\mid x_{N})(y^{*})\mid(x^{*},y^{*})\in\partial\varphi(x_{0},x_{N})\},\]
where $G_{\{\tilde{x}_{i}\}_{i=0}^{N}}:\mathbb{R}^{n}\rightrightarrows\mathbb{R}^{n}$
and $\mathcal{X}_{N}$ are defined as in \eqref{eq:D-R-N-first-line}
and \eqref{eq:X-N}. If in addition all the $F_{k,N}(\cdot)$ are
all graph convex and $\varphi(\cdot,\cdot)$ is convex, we have \begin{eqnarray*}
\partial f(x_{0}) & = & \{x^{*}+D^{*}G_{\{\tilde{x}_{i}\}_{i=0}^{N}}(x_{0}\mid x_{N})(y^{*})\mid(x^{*},y^{*})\in\partial\varphi(x_{0},x_{N})\}\\
 &  & \mbox{for any }\{\tilde{x}_{i}\}_{i=0}^{N}\in\mathcal{X}_{N}\mbox{ s.t. }f(x_{0})=\varphi(x_{0},x_{N}).\end{eqnarray*}

In particular, a necessary condition for the optimality of the path
$\{\tilde{x}_{i}\}_{i=0}^{N}\in\mathcal{X}_{N}$ is the existence
of $\{p_{i}\}_{i=0}^{N}$ such that 
\begin{enumerate}
\item $(-p_{0},p_{N})\in\partial\varphi(\tilde{x}_{0},\tilde{x}_{N})$,
and
\item \textup{\emph{$\frac{p_{k}-p_{k-1}}{\Delta t}\in-D^{*}F_{k-1,N}\left(\tilde{x}_{k-1}\mid\frac{1}{\Delta t}(\tilde{x}_{k}-\tilde{x}_{k-1})\right)(p_{k})\mbox{ for all }k\in\overline{1,N}.$}}
\end{enumerate}
\end{thm}
\begin{proof}
Apply Theorem \ref{thm:discrete-coderv} and Lemma \ref{lem:gen-opt}.
\end{proof}
Condition (1) in Theorem \ref{thm:nec-condn-discrete} is the discrete
analogue of the Transversality Condition (TC), while condition (2)
is the analogue of the Euler-Lagrange condition (EL). 

Finally, we make a remark on the Weierstrass-Pontryagin Maximum Principle
(WP). Before we do so, we recall that for $F:[0,T]\times\mathbb{R}^{n}\rightrightarrows\mathbb{R}^{n}$,
the reachable map of the \emph{relaxed }differential inclusion is
defined by \begin{eqnarray*}
R_{\scriptsize{\overline{\co}F}}(x_{0}) & := & \{y:\exists x(\cdot)\in AC([0,T],\mathbb{R}^{n})\mbox{ s.t. }\\
 &  & \qquad x^{\prime}(t)\in\overline{\co}F\big(t,x(t)\big)\mbox{ for }t\in[0,T]\mbox{ a.e.},\\
 &  & \qquad x(0)=x_{0}\mbox{ and }x(T)=y\}.\end{eqnarray*}
It is well known that under mild conditions, we have $\cl\, R(x)=R_{\scriptsize{\overline{\co}}F}(x)$
for all $x\in\mathbb{R}^{n}$. 
\begin{rem}
(Discrete analogue of the Weierstrass-Pontryagin Maximum Principle)
Recall the chain rule for set-valued maps $F:\mathbb{R}^{n}\rightrightarrows\mathbb{R}^{n}$
and $G:\mathbb{R}^{n}\rightrightarrows\mathbb{R}^{n}$ as presented
in Theorem \ref{thm:Coderivative-chain-rule}. If the conclusion of
the chain rule had been that for all $r\in\mathbb{R}^{n}$, \begin{eqnarray}
 &  & D^{*}(F\circ G)(\bar{x}\mid\bar{z})(r)\label{eq:co-chain-rule-3}\\
 & \subset & \bigcup_{\bar{y}\in F^{-1}(\bar{z})\cap G(\bar{x})}\{\overline{\co}D^{*}G(\bar{x}\mid\bar{y})(q)|q\in\overline{\co}D^{*}F(\bar{x}\mid\bar{y})(r),\left\langle q,\bar{y}-y^{\prime}\right\rangle \leq0\mbox{ for all }y^{\prime}\in G(\bar{x})\},\nonumber \end{eqnarray}
then we can repeatedly apply this chain rule like in Theorem \ref{thm:discrete-coderv}
so that under the conditions of Theorem \ref{thm:discrete-coderv},
$p_{0}\in D^{*}R_{N}(x_{0}|x_{N})(p_{N})$ implies that there are
$\{\tilde{x}_{i}\}_{i=0}^{N}$ and $\{\tilde{p}_{i}\}_{i=0}^{N}$
such that \begin{subequations} \begin{eqnarray}
 &  & \tilde{x}_{0}=x_{0},\,\tilde{x}_{N}=x_{N},\,\tilde{p}_{0}=p_{0},\,\tilde{p}_{N}=p_{N},\label{eq:WP-1}\\
 &  & \frac{\tilde{p}_{k}-\tilde{p}_{k-1}}{\Delta t}\in-D^{*}F_{k-1,N}\left(\tilde{x}_{k-1}\mid\frac{1}{\Delta t}(\tilde{x}_{k}-\tilde{x}_{k-1})\right)(\tilde{p}_{k})\label{eq:WP-2}\\
 &  & \qquad\qquad\mbox{ for all }k\in\overline{1,N},\nonumber \\
 & \mbox{ and } & \left\langle -\tilde{p}_{k},v-\frac{1}{\Delta t}(\tilde{x}_{k}-\tilde{x}_{k-1})\right\rangle \leq0\label{eq:WP-3}\\
 &  & \qquad\qquad\mbox{ for all }v\in F_{k-1,N}(\tilde{x}_{k-1})\mbox{ and }k\in\overline{1,N}.\nonumber \end{eqnarray}
\end{subequations}Such a formula would be appealing because \eqref{eq:WP-2}
corresponds to the Euler-Lagrange Condition (EL) and \eqref{eq:WP-3}
corresponds to the Weierstrass-Pontryagin Maximum Principle (WP).
However, \eqref{eq:co-chain-rule-3} is not true in general. Consider
the maps $G:\mathbb{R}\rightrightarrows\mathbb{R}$ and $f:\mathbb{R}\to\mathbb{R}$
defined by \begin{eqnarray*}
G(x) & := & [x+1,x+2]\cup[x-2,x-1],\\
\mbox{ and }f(x) & := & -|x-0.5|.\end{eqnarray*}
Then $f\circ G(0)=[-2.5,-0.5]$, and $f^{-1}(-0.5)\cap G(0)=\{1\}$.
We can calculate that \begin{eqnarray*}
 &  & D^{*}G(0\mid1)(1)=\{1\},\\
 &  & D^{*}f(1\mid-0.5)(-1)=\{1\},\\
 & \mbox{and } & D^{*}(f\circ G)(0|-0.5)(-1)=\{1\}.\end{eqnarray*}
However, since we do not have$\left\langle 1,1-v\right\rangle \leq0$
for all $v\in G(0)=[1,2]\cup[-2,-1]$, the right hand side of \eqref{eq:co-chain-rule-3}
is empty, showing us that \eqref{eq:co-chain-rule-3} cannot be true.

In the case where $F_{k,N}:\mathbb{R}^{n}\rightrightarrows\mathbb{R}^{n}$
are convex-valued (so that we are considering the relaxed differential
inclusion), it is an easy exercise that provided $F_{k,N}$ is continuous,
then \eqref{eq:WP-2} is equivalent to \[
\left(-\frac{\tilde{p}_{k}-\tilde{p}_{k-1}}{\Delta t},-\tilde{p}_{k}\right)\in N_{\scriptsize{\gph(F_{k-1,N})}}\left(\tilde{x}_{k-1}\mid\frac{1}{\Delta t}(\tilde{x}_{k}-\tilde{x}_{k-1})\right).\]
In addition to the fact that $F_{k-1,N}(\tilde{x}_{k-1})$ is convex,
\eqref{eq:WP-3} follows easily.
\end{rem}

\section{\label{sec:cts-incl}Subdifferential analysis of differential inclusions}

In this section, we make use of the work in Section \ref{sec:discrete-incl}
to calculate estimates of the convexified limiting coderivative of
the (continuous) reachable map, and explain how this new formula gives
a new way to interpret the Euler-Lagrange and Transversality conditions. 

We first simplify the notation. Define $\mathcal{F}(x,y)$ to be the
set of feasible paths with end points $\bar{x}$ and $\bar{y}$, i.e.,
\begin{eqnarray}
\mathcal{F}(x,y) & := & \{x(\cdot)|x(\cdot)\in AC([0,T],\mathbb{R}^{n}),x(0)=x,x(T)=y,\label{eq:F_x_y}\\
 &  & \qquad\mbox{and }x^{\prime}(t)\in F(t,x(t))\mbox{ for }t\in[0,T]\mbox{ a.e.}\}.\nonumber \end{eqnarray}
Define $\Pi:\mathbb{R}^{n}\times\mathbb{R}^{n}\times\mathbb{R}^{n}\rightrightarrows\mathbb{R}^{n}$
by \begin{eqnarray*}
\Pi(x,y,v) & := & \{u\mid\exists x(\cdot)\in\mathcal{F}(x,y),p(\cdot)\in AC([0,T],\mathbb{R}^{n})\\
 &  & \qquad\mbox{s.t. }p(0)=u,\, p(T)=v\mbox{ and }\\
 &  & \qquad p^{\prime}(t)\in-\overline{\co}D_{x}^{*}F\big(t,x(t)\mid x^{\prime}(t)\big)\big(p(t)\big)\mbox{ for }t\in[0,T]\mbox{ a.e.}\}.\end{eqnarray*}
Here, $\overline{\co}D_{x}^{*}F(t,x(t)\mid x^{\prime}(t)):\mathbb{R}^{n}\rightrightarrows\mathbb{R}^{n}$
is to be understood as $\overline{\co}D^{*}F_{t}(x(t)\mid x^{\prime}(t)):\mathbb{R}^{n}\rightrightarrows\mathbb{R}^{n}$,
where $F_{t}(\cdot)=F(t,\cdot)$. Corresponding to $\Pi(x,y,v)$ is
its discretized version: \begin{eqnarray*}
 &  & \Pi_{N}(x,y,v)\\
 & := & \bigg\{u\mid\exists\{x_{i}\}_{i=0}^{N},\{p_{i}\}_{i=0}^{N}\mbox{ s.t. }x_{0}=x,\, x_{N}=y,\, p_{0}=u,\, p_{N}=v,\\
 &  & \qquad\frac{1}{\Delta t}(x_{k}-x_{k-1})\in F\big((k-1)(\Delta t),x_{k-1}\big)\mbox{ for all }k\in\overline{1,N},\mbox{ and }\\
 &  & \qquad\frac{1}{\Delta t}(p_{k}-p_{k-1})\in-D^{*}F_{k-1,N}\left(x_{k-1}\mid\frac{1}{\Delta t}(x_{k}-x_{k-1})\right)(p_{k})\mbox{ for all }k\in\overline{1,N}\bigg\}.\end{eqnarray*}
We make the following conjecture. 
\begin{conjecture}
\label{con:cts-coderv}(Upper estimate of discretized coderivative
of reachable map) For the reachable map $R:\mathbb{R}^{n}\rightrightarrows\mathbb{R}^{n}$
defined in \eqref{eq:R-diff-reachable}, the convexified coderivative
$\overline{\co}D^{*}R(\bar{x}\mid\bar{y})$ satisfies \begin{eqnarray}
D^{*}R(\bar{x}\mid\bar{y})(v) & \subset & \{u:\exists x(\cdot)\in\mathcal{F}(\bar{x},\bar{y}),p(\cdot)\in AC([0,T],\mathbb{R}^{n})\mbox{ s.t. }\label{eq:R-coderv-conj}\\
 &  & \qquad p^{\prime}(t)\in-\overline{\co}D_{x}^{*}F\big(t,x(t)\mid x^{\prime}(t)\big)\big(p(t)\big),\nonumber \\
 &  & \qquad p_{0}=u\mbox{ and }p_{T}=v\}\mbox{ for all }v\in\mathbb{R}^{n}\}.\nonumber \end{eqnarray}
\end{conjecture}
\begin{rem}
(Consequence of Conjecture \ref{con:cts-coderv}) Consider the problem
\begin{eqnarray*}
 & \underset{\scriptsize{(x,y)}}{\min} & \varphi(x,y)\\
 & \mbox{s.t.} & (x,y)\in\gph(R).\end{eqnarray*}
Recall the discussion in Remark \ref{rem:alternate-marginal}. Provided
\eqref{eq:R-coderv-conj} holds, if the point $(\bar{x},\bar{y})$
is optimal, then $0\in\partial\varphi(\bar{x},\bar{y})+N_{\scriptsize\gph(R)}(\bar{x},\bar{y})$.
We have \begin{eqnarray*}
\partial\varphi(\bar{x},\bar{y})+N_{\scriptsize\gph(R)}(\bar{x},\bar{y}) & = & \partial\varphi(\bar{x},\bar{y})+L\,\gph\big(D^{*}R(\bar{x},\bar{y})\big)\\
 & \subset & \partial\varphi(\bar{x},\bar{y})+L\,\gph\big(\overline{\co}D^{*}R(\bar{x}\mid\bar{y})\big),\end{eqnarray*}
where $L:\mathbb{R}^{n}\times\mathbb{R}^{n}\to\mathbb{R}^{n}\times\mathbb{R}^{n}$
is the linear map represented by the matrix $\left({0\atop -I}{I\atop 0}\right)$.
Unrolling the definition of $D^{*}R(\bar{x}\mid\bar{y})$ gives the
following optimality condition: If $(\bar{x},\bar{y})$ is optimal,
then there are paths $x(\cdot),p(\cdot)\in AC([0,T],\mathbb{R}^{n})$
such that $x(\cdot)$ is feasible for the differential inclusion,
$x(0)=\bar{x}$, $x(T)=\bar{y}$, and satisfies the Transversality
condition (TC) and the Euler-Lagrange condition (EL).
\end{rem}
We will prove the following weaker result instead: \begin{eqnarray}
\overline{\co}D^{*}R(\bar{x}\mid\bar{y})(v) & \subset & \overline{\co}\{u:\exists x(\cdot)\in\mathcal{F}(\bar{x},\bar{y}),p(\cdot)\in AC([0,T],\mathbb{R}^{n})\mbox{ s.t. }\label{eq:R-conv-coderv}\\
 &  & \qquad p^{\prime}(t)\in-\overline{\co}D_{x}^{*}F\big(t,x(t)\mid x^{\prime}(t)\big)\big(p(t)\big),\nonumber \\
 &  & \qquad p_{0}=u\mbox{ and }p_{T}=v\}\mbox{ for all }v\in\mathbb{R}^{n}.\nonumber \end{eqnarray}
Our strategy is to prove the following three inclusions:\begin{subequations}
\begin{eqnarray}
\overline{\co}D^{*}R(\bar{x}\mid\bar{y}) & \subset & \bigcap_{{N\in\mathbb{N}\atop \delta>0}}\overline{\co}\bigcup_{i>N}\bigcup_{{x\in\mathbb{B}_{\delta}(\bar{x})\atop y\in\mathbb{B}_{\delta}(\bar{y})}}D^{*}R_{i}(x\mid y),\label{eq:conj-step1}\\
\bigcap_{{N\in\mathbb{N}\atop \delta>0}}\overline{\co}\bigcup_{i>N}\bigcup_{{x\in\mathbb{B}_{\delta}(\bar{x})\atop y\in\mathbb{B}_{\delta}(\bar{y})}}D^{*}R_{i}(x\mid y)(v) & \subset & \bigcap_{{N\in\mathbb{N}\atop \delta>0}}\overline{\co}\bigcup_{i>N}\bigcup_{{x\in\mathbb{B}_{\delta}(\bar{x})\atop y\in\mathbb{B}_{\delta}(\bar{y})}}\Pi_{i}(x,y,v),\label{eq:conj-step2}\\
\mbox{ and }\bigcap_{{N\in\mathbb{N}\atop \delta>0}}\overline{\co}\bigcup_{i>N}\bigcup_{{x\in\mathbb{B}_{\delta}(\bar{x})\atop y\in\mathbb{B}_{\delta}(\bar{y})}}\Pi_{i}(x,y,v) & \subset & \overline{\co}\Pi(\bar{x},\bar{y},v),\label{eq:conj-step3}\end{eqnarray}
\end{subequations}where \eqref{eq:conj-step2} and \eqref{eq:conj-step3}
hold for all $v\in\mathbb{R}^{n}$. Conditions for $D^{*}R_{i}(x\mid y)(v)\subset\Pi_{p,i}(x,y,v)$,
which addresses \eqref{eq:conj-step2}, were discussed in Theorem
\ref{thm:discrete-coderv}. The same steps used to prove that \eqref{eq:conj-step2}
and \eqref{eq:conj-step3} hold for all $v\in\mathbb{R}^{n}$ yield
the following stronger statements: For all $v\in\mathbb{R}^{n}$,
\begin{eqnarray*}
\bigcap_{{N\in\mathbb{N}\atop \delta>0}}\cl\bigcup_{i>N}\bigcup_{{x\in\mathbb{B}_{\delta}(\bar{x})\atop y\in\mathbb{B}_{\delta}(\bar{y})}}D^{*}R_{i}(x\mid y)(v) & \subset & \bigcap_{{N\in\mathbb{N}\atop \delta>0}}\cl\bigcup_{i>N}\bigcup_{{x\in\mathbb{B}_{\delta}(\bar{x})\atop y\in\mathbb{B}_{\delta}(\bar{y})}}\Pi_{i}(x,y,v),\\
\mbox{ and }\bigcap_{{N\in\mathbb{N}\atop \delta>0}}\cl\bigcup_{i>N}\bigcup_{{x\in\mathbb{B}_{\delta}(\bar{x})\atop y\in\mathbb{B}_{\delta}(\bar{y})}}\Pi_{i}(x,y,v) & \subset & \Pi(\bar{x},\bar{y},v).\end{eqnarray*}
Notice that if \eqref{eq:conj-step1} were strengthened to be \[
D^{*}R(\bar{x}\mid\bar{y})\subset\bigcap_{{N\in\mathbb{N}\atop \delta>0}}\cl\bigcup_{i>N}\bigcup_{{x\in\mathbb{B}_{\delta}(\bar{x})\atop y\in\mathbb{B}_{\delta}(\bar{y})}}D^{*}R_{i}(x\mid y)\]
instead, then piecing the last three formulas together gives \eqref{eq:R-coderv-conj}.
We continue with some lemmas.
\begin{lem}
\label{lem:Coderv-ard-x-y}(Coderivatives around $(\bar{x},\bar{y})$)
Let $\delta>0$, and $S:\mathbb{R}^{n}\rightrightarrows\mathbb{R}^{m}$
be a closed set-valued map. Suppose $H:\mathbb{R}^{n}\rightrightarrows\mathbb{R}^{m}$
is a prefan such that \[
H\in\mathcal{H}\left(\overline{\co}\bigcup_{{x\in\mathbb{B}_{\delta}(\bar{x})\atop y\in\mathbb{B}_{\delta}(\bar{y})}}D^{*}S(x\mid y)\right).\]
Let $\delta^{\prime}:=\min(\delta,\frac{\delta}{5\|H\|^{+}})$. Then
\[
S(x^{\prime})\cap\mathbb{B}_{\delta/2}(\bar{y})\subset S(x^{\prime\prime})+H(x^{\prime}-x^{\prime\prime})\mbox{ for all }x^{\prime},x^{\prime\prime}\in\mathbb{B}_{\delta^{\prime}}(\bar{x}).\]
\end{lem}
\begin{proof}
For any $x\in\mathbb{B}_{\delta}(\bar{x})$ and $y\in\mathbb{B}_{\delta}(\bar{y})$,
we have $H\in\mathcal{H}(D^{*}S(x\mid y))$. Choose any $\theta>0$.
There exists some $\epsilon_{x,y,\theta}>0$ such that \[
S(x^{\prime})\cap\mathbb{B}_{\epsilon_{x,y,\theta}}(y)\subset S(x^{\prime\prime})+(H+\theta)(x^{\prime}-x^{\prime\prime})\mbox{ for all }x^{\prime},x^{\prime\prime}\in\mathbb{B}_{\epsilon_{x,y,\theta}}(x).\]
For each $x\in\mathbb{B}_{\delta}(\bar{x})$, we can find a finite
number of elements in $\mathbb{B}_{\delta}(\bar{y})$, say $\{y_{j}\}_{j=1}^{J}$,
such that $\mathbb{B}_{\delta}(\bar{y})\subset\cup_{j=1}^{J}\mathbb{B}_{\epsilon_{x,y,\theta}}(y)$.
Letting $\epsilon_{x,\theta}:=\min_{j\in\overline{1,J}}(\epsilon_{x,y_{j},\theta})$,
we have \[
S(x^{\prime})\cap\mathbb{B}_{\delta}(\bar{y})\subset S(x^{\prime\prime})+(H+\theta)(x^{\prime}-x^{\prime\prime})\mbox{ for all }x^{\prime},x^{\prime\prime}\in\mathbb{B}_{\epsilon_{x,\theta}}(x).\]
For any line segment $[x^{\prime},x^{\prime\prime}]$ in $\mathbb{B}_{\delta}(\bar{x})$,
we can find finitely many $x$ in $\mathbb{B}_{\delta}(\bar{x})$,
say $\{x_{k}\}_{k=1}^{K}$ such that $[x^{\prime},x^{\prime\prime}]\subset\cup_{k=1}^{K}\mathbb{B}_{\epsilon_{x_{k},\theta}}(x_{k})$.
We can break up the line segment $[x^{\prime},x^{\prime\prime}]$
to a union of line segments $\cup_{j=1}^{J-1}[\tilde{x}_{j},\tilde{x}_{j+1}]$
so that $\{\tilde{x}_{j}\}_{j=1}^{J}$ line up in that order, each
$[\tilde{x}_{j},\tilde{x}_{j+1}]$ is inside some $\mathbb{B}_{\epsilon_{x_{k},\theta}}(x_{k})$,
$\tilde{x}_{1}=x^{\prime}$ and $\tilde{x}_{J}=x^{\prime\prime}$.
Then\begin{eqnarray*}
S(\tilde{x}_{j})\cap\mathbb{B}_{\delta}(\bar{y}) & \subset & S(\tilde{x}_{j+1})+(H+\theta)(\tilde{x}_{j}-\tilde{x}_{j+1})\\
\Rightarrow[S(\tilde{x}_{j})\cap\mathbb{B}_{\delta}(\bar{y})]+(H+\theta)(\tilde{x}_{1}-\tilde{x}_{j}) & \subset & S(\tilde{x}_{j+1})+(H+\theta)(\tilde{x}_{1}-\tilde{x}_{j+1}).\end{eqnarray*}
We write $\kappa=\|H\|^{+}$ to simplify notation. This gives\begin{eqnarray*}
[S(\tilde{x}_{j})+(H+\theta)(\tilde{x}_{1}-\tilde{x}_{j})]\cap\mathbb{B}_{\delta-(\kappa+\theta)|x^{\prime}-x^{\prime\prime}|}(\bar{y}) & \subset & [S(\tilde{x}_{j})\cap\mathbb{B}_{\delta}(\bar{y})]+(H+\theta)(\tilde{x}_{1}-\tilde{x}_{j})\\
 & \subset & S(\tilde{x}_{j+1})+(H+\theta)(\tilde{x}_{1}-\tilde{x}_{j+1}),\end{eqnarray*}
which implies \begin{eqnarray}
 &  & [S(\tilde{x}_{j})+(H+\theta)(\tilde{x}_{1}-\tilde{x}_{j})]\cap\mathbb{B}_{\delta-(\kappa+\theta)|x^{\prime}-x^{\prime\prime}|}(\bar{y})\label{eq:gen-diff-final}\\
 & \subset & [S(\tilde{x}_{j+1})+(H+\theta)(\tilde{x}_{1}-\tilde{x}_{j+1})]\cap\mathbb{B}_{\delta-(\kappa+\theta)|x^{\prime}-x^{\prime\prime}|}(\bar{y}).\nonumber \end{eqnarray}
Consider the case where $\theta<\kappa/4$ so that $4(\kappa+\theta)<5\kappa$.
If $x^{\prime},x^{\prime\prime}\in\mathbb{B}_{\delta^{\prime}}(\bar{x})$,
where $\delta^{\prime}=\min(\delta,\frac{\delta}{5\kappa})$, then
\[
(\kappa+\theta)|x^{\prime}-x^{\prime\prime}|\leq\frac{5}{4}\kappa\frac{2\delta}{5\kappa}\leq\frac{\delta}{2}.\]
Recalling that $\tilde{x}_{1}=x^{\prime}$ and $\tilde{x}_{J}=x^{\prime\prime}$
and applying \eqref{eq:gen-diff-final} repeatedly, we have \[
S(x^{\prime})\cap\mathbb{B}_{\delta/2}(\bar{y})\subset S(x^{\prime\prime})+(H+\theta)(x^{\prime}-x^{\prime\prime}).\]
The above holds for all $x^{\prime},x^{\prime\prime}\in\mathbb{B}_{\delta^{\prime}}(\bar{x})$
and for all $\theta>0$, and hence for $\theta=0$, giving us the
conclusion we need.
\end{proof}
This result gives a handle on the left hand bound.
\begin{lem}
\label{lem:On-conj-1}(On \eqref{eq:conj-step1}) Let $S:\mathbb{R}^{n}\rightrightarrows\mathbb{R}^{m}$
be a closed set-valued map. Suppose $\{S_{i}(\cdot)\}_{i=1}^{\infty}$,
where $S_{i}:\mathbb{R}^{n}\rightrightarrows\mathbb{R}^{m}$, are
closed set-valued maps such that for any $\epsilon>0$ and $x\in\mathbb{R}^{n}$,
there is some $I$ such that \begin{equation}
\mathbf{d}(S(x),S_{i}(x))<\epsilon\mbox{ for all }i>I.\label{eq:Hausdorff-bdd}\end{equation}
Then for any $\delta>0$ and positive integer $N$, we have \[
\overline{\co}D^{*}S(\bar{x}\mid\bar{y})\subset\overline{\co}\bigcup_{i>N}\bigcup_{{x\in\mathbb{B}_{\delta}(\bar{x})\atop y\in\mathbb{B}_{\delta}(\bar{y})}}D^{*}S_{i}(x\mid y).\]
\end{lem}
\begin{proof}
By Lemma \ref{lem:conv-coderv-gen-derv}, we can prove that the following
holds for all $\delta>0$ and positive integers $N$ instead: \[
\mathcal{H}\big(\overline{\co}D^{*}S(\bar{x}\mid\bar{y})\big)\supset\mathcal{H}\left(\overline{\co}\bigcup_{i>N}\bigcup_{{x\in\mathbb{B}_{\delta}(\bar{x})\atop y\in\mathbb{B}_{\delta}(\bar{y})}}D^{*}S_{i}(x\mid y)\right).\]
Suppose $H:\mathbb{R}^{n}\rightrightarrows\mathbb{R}^{m}$ is a prefan
in the RHS. Then  for any $i>N$ and $\delta>0$, \[
H\in\mathcal{H}\left(\overline{\co}\bigcup_{{x\in\mathbb{B}_{\delta}(\bar{x})\atop y\in\mathbb{B}_{\delta}(\bar{y})}}D^{*}S_{i}(x\mid y)\right).\]
By Lemma \ref{lem:Coderv-ard-x-y}, if $\delta^{\prime}=\min(\delta,\frac{\delta}{5\|H\|^{+}})$
, then \[
S_{i}(x^{\prime})\cap\mathbb{B}_{\delta/2}(\bar{y})\subset S_{i}(x^{\prime\prime})+H(x^{\prime}-x^{\prime\prime})\mbox{ for all }x^{\prime},x^{\prime\prime}\in\mathbb{B}_{\delta^{\prime}}(\bar{x}).\]
For all $x^{\prime},x^{\prime\prime}\in\mathbb{B}_{\delta^{\prime}}(\bar{x})$
and $\epsilon>0$, we can find $i$ large enough so that \begin{eqnarray*}
S(x^{\prime})\cap\mathbb{B}_{\delta/2}(\bar{y}) & \subset & [S_{i}(x^{\prime})+\epsilon\mathbb{B}]\cap\mathbb{B}_{\delta/2}(\bar{y})\\
 & \subset & S_{i}(x^{\prime\prime})+H(x^{\prime}-x^{\prime\prime})+\epsilon\mathbb{B}\\
 & \subset & S(x^{\prime\prime})+H(x^{\prime}-x^{\prime\prime})+2\epsilon\mathbb{B}.\end{eqnarray*}
The above holds for all $\epsilon>0$, and we have \[
S(x^{\prime})\cap\mathbb{B}_{\delta/2}(\bar{y})\subset S(x^{\prime\prime})+H(x^{\prime}-x^{\prime\prime})\mbox{ for all }x^{\prime},x^{\prime\prime}\in\mathbb{B}_{\delta^{\prime}}(\bar{x}).\]
This implies that $H\in\mathcal{H}\big(\overline{\co}D^{*}S(\bar{x}\mid\bar{y})\big)$
as needed.\end{proof}
\begin{rem}
(On formula \eqref{eq:Hausdorff-bdd}) We note that conditions for
$\mathbf{d}(S(x),S_{i}(x))<\epsilon$ were given in \cite{DL92},
and in particular, conditions for $S(x)\subset S_{i}(x)+\epsilon\mathbb{B}$
were given in \cite[Theorem 6.4]{Mor06} for example.
\end{rem}
Note that Theorem \ref{thm:discrete-coderv} says that $D^{*}S_{i}(x\mid y)(v)\subset\Pi_{i}(x,y,v)$.
To find suitable conditions for \eqref{eq:conj-step3}, we need the
following result.
\begin{lem}
\label{lem:nested-sets}(Convexification of intersection of nested
sets) Suppose $\{A_{i}\}_{i=1}^{\infty}\subset\mathbb{R}^{n}$ are
nested compact sets such that $A_{i+1}\subset A_{i}$. Then $\overline{\co}\cap_{i}A_{i}=\cap_{i}\overline{\co}A_{i}$. \end{lem}
\begin{proof}
Suppose $x$ is in the LHS. Then $x\in\overline{\co}A_{i}$ for all
$i$, so $x\in\cap_{i}\overline{\co}A_{i}$, establishing $\overline{\co}\cap_{i}A_{i}\subset\cap_{i}\overline{\co}A_{i}$.

Next, suppose $x$ is in the RHS. Then $x\in\overline{\co}A_{i}$
for all $i$. Consider any $v\in\mathbb{R}^{n}\backslash\{0\}$. Since
$x\in\overline{\co}A_{i}$, we have $v^{T}x\leq\sup_{a\in A_{i}}v^{T}a$.
By the compactness of $A_{i}$, let $\bar{a}_{i}$ be such that $v^{T}\bar{a}_{i}=\sup_{a\in A_{i}}v^{T}a$.
Since $\cap_{j}A_{j}\subset A_{i}$ for all $i$, it is clear that
$\sup_{a\in\cap_{j}A_{j}}v^{T}a\leq\sup_{a\in A_{i}}v^{T}a$ for all
$i$, so $\sup_{a\in\cap_{j}A_{j}}v^{T}a\leq\inf_{i}\sup_{a\in A_{i}}v^{T}a$.
By the compactness of $A_{i}$, the limit $\bar{a}=\lim_{j\to\infty}\bar{a}_{j}$
exists and lies in $\cap_{j}A_{j}$. This shows that \begin{eqnarray*}
\inf_{i}\sup_{a\in A_{i}}v^{T}a & = & \inf_{i}v^{T}\bar{a}_{i}\\
 & = & v^{T}\bar{a}\\
 & \leq & \sup_{a\in\cap_{j}A_{j}}v^{T}a.\end{eqnarray*}
Then $v^{T}x\leq\sup_{a\in\cap_{j}A_{j}}v^{T}a$, which holds for
all $v$. Thus we have $x\in\overline{\co}\cap_{i}A_{i}$, so $\overline{\co}\cap_{i}A_{i}=\cap_{i}\overline{\co}A_{i}$
as needed.
\end{proof}
Here is a lemma useful for proving our next result. We take our result
from \cite[Lemma 4.4]{Smirnov02}.
\begin{lem}
\label{lem:Smirnov-L4.4}(Continuous solutions from discrete solutions)
Assume that a set-valued map $F:[0,T]\times\mathbb{R}^{n}\times\mathbb{R}^{m}\to\mathbb{R}^{n}$
has closed convex values. Let the set-valued map $(x,y)\mapsto F(t,x,y)$
be upper semicontinuous for almost all $t\in[0,T]$, and let $F(t,x,y)\subset b(t)\mathbb{B}$
for all $(t,x,y)\in[0,T]\times\mathbb{R}^{n}\times\mathbb{R}^{m}$,
where $b(\cdot)\in L_{1}([0,T],\mathbb{R})$. Assume that functions
$x_{k}(\cdot)\in AC([0,T],\mathbb{R}^{n})$, $k=0,1,\dots,$ satisfy
\[
x_{k}^{\prime}(t)\in\overline{\co}F\big(t,x_{k}(t),\eta_{k}(t)\mathbb{B}_{m}\big)+\eta_{k}(t)\mathbb{B}_{n},\]
where $\eta_{k}\geq0$, $\lim_{k\to\infty}\eta_{k}(t)=0$ almost everywhere,
and $|\eta_{k}(t)|\leq\eta(t)$, $k=1,2,\dots$, $\eta(\cdot)\in L_{1}([0,T],\mathbb{R})$.
Then the functions $x_{k}(\cdot)$ are equicontinuous on $[0,T]$;
and if a subsequence $x_{k_{p}}(\cdot)$ uniformly converges to a
function $x(\cdot)$, then $x(\cdot)$ is a solution of the differential
inclusion \[
x^{\prime}(t)\in F\big(t,x(t),0\big)\mbox{ for }t\in[0,T]\mbox{ a.e.}\]

\end{lem}
Before we state our main result, we describe in detail the paths produced
by discrete approximations in the remark below.
\begin{rem}
\label{rem:Disc-approx}(Discrete approximations) For $\{x_{i,j}\}_{{0\leq j\leq i\atop 1\leq i\leq\infty}}$
and $\{p_{i,j}\}_{{0\leq j\leq i\atop 1\leq i\leq\infty}}$, let $\Delta t=T/i$,
and construct the following path $x_{i}:[0,T]\to\mathbb{R}^{n}$ defined
by \begin{eqnarray*}
x_{i}(t) & = & \frac{t-j\Delta t}{\Delta t}x_{i,j+1}+\frac{(j+1)\Delta t-t}{\Delta t}x_{i,j}\mbox{ whenever }t\in[j\Delta t,(j+1)\Delta t],\\
\mbox{ and }p_{i}(t) & = & \frac{t-j\Delta t}{\Delta t}p_{i,j+1}+\frac{(j+1)\Delta t-t}{\Delta t}p_{i,j}\mbox{ whenever }t\in[j\Delta t,(j+1)\Delta t].\end{eqnarray*}
It is clear that $x_{i}(\cdot)$ and $p_{i}(\cdot)$ are piecewise
differentiable at all points other than integer multiples of $\Delta t$,
and the derivatives satisfy \begin{subequations} \begin{eqnarray}
x_{i}^{\prime}(t) & = & \frac{1}{\Delta t}(x_{i,j+1}-x_{i,j})\mbox{ whenever }t\in\big(j\Delta t,(j+1)\Delta t\big),\label{eq:x-i-prime}\\
\mbox{ and }p_{i}^{\prime}(t) & = & \frac{1}{\Delta t}(p_{i,j+1}-p_{i,j})\mbox{ whenever }t\in\big(j\Delta t,(j+1)\Delta t\big).\label{eq:p-i-prime}\end{eqnarray}

\end{subequations}We also need the following condition for Lemma
\ref{lem:3rd-containment}, which was one of the conclusions in Theorem
\ref{thm:nec-condn-discrete}: \emph{\begin{equation}
\frac{p_{k}-p_{k-1}}{\Delta t}\in-D^{*}F_{k-1,N}\left(x_{k-1}\mid\frac{1}{\Delta t}(x_{k}-x_{k-1})\right)(p_{k})\mbox{ for all }k\in\overline{1,N}.\label{eq:disc-diff-incl-condn}\end{equation}
}
\end{rem}
We now prove our result on \eqref{eq:conj-step3}. Note that \eqref{eq:conj-step3}
represents a closedness property, and we shall show that Lemma \ref{lem:Smirnov-L4.4}
provides some reasonable conditions for \eqref{eq:conj-step3} to
hold.
\begin{lem}
\label{lem:3rd-containment}(On \eqref{eq:conj-step3}) Suppose $F:[0,T]\times\mathbb{R}^{n}\rightrightarrows\mathbb{R}^{n}$
is a convex-valued osc function. Assume further that there is some
$b(\cdot)\in L_{1}([0,T],\mathbb{R}^{n})$ such that  $\overline{\co}D_{x}^{*}F(t,x\mid y)(p)\subset b(t)\|p\|\mathbb{B}$
for all $(t,x,y,p)\in[0,T]\times\mathbb{R}^{n}\times\mathbb{R}^{n}\times\mathbb{R}^{n}$.
Suppose also that the following assumption holds:
\begin{enumerate}
\item Whenever $\{x_{i}(\cdot)\}_{i=1}^{\infty}$ and $\{p_{i}(\cdot)\}_{i=1}^{\infty}$,
constructed based on discrete approximations $\{x_{i,j}\}_{{0\leq j\leq i\atop 1\leq i\leq\infty}}$
and $\{p_{i,j}\}_{{0\leq j\leq i\atop 1\leq i\leq\infty}}$ satisfying
\eqref{eq:disc-diff-incl-condn} as described in Remark \ref{rem:Disc-approx},
satisfies $x_{i}(0)\to\bar{x}$ and $x_{i}(T)\to\bar{y}$ as $N\to\infty$,
then there exists a subsequence of $\{x_{i}(\cdot)\}_{i=1}^{\infty}$,
say $\{x_{i_{k}}(\cdot)\}_{k=1}^{\infty}$, and $x(\cdot)\in AC([0,T],\mathbb{R}^{n})$
such that

\begin{itemize}
\item $x_{i_{k}}(\cdot)$ converges uniformly to $x(\cdot)$, 
\item $x_{i_{k}}^{\prime}(\cdot)$ converges pointwise almost everywhere
to $x^{\prime}(\cdot)$,
\item $x(\cdot)$ satisfies the differential inclusion \begin{eqnarray*}
 &  & x^{\prime}(t)\in F\big(t,x(t)\big)\mbox{ a.e.},\\
 &  & x(0)=\bar{x}\mbox{ and }x(T)=\bar{y},\end{eqnarray*}

\item and $\{p_{i_{k}}(\cdot)\}_{k=1}^{\infty}$ converges uniformly to
some $p(\cdot)$. 
\end{itemize}
\end{enumerate}
Then we have

\begin{eqnarray}
\bigcap_{N\in\mathbb{N},\delta>0}\cl\bigcup_{i>N}\bigcup_{{x\in\mathbb{B}_{\delta}(\bar{x})\atop y\in\mathbb{B}_{\delta}(\bar{y})}}\Pi_{i}(x,y,v) & \subset & \Pi(\bar{x},\bar{y},v).\label{eq:nonconvex-target}\\
\mbox{ and }\bigcap_{N\in\mathbb{N},\delta>0}\overline{\co}\bigcup_{i>N}\bigcup_{{x\in\mathbb{B}_{\delta}(\bar{x})\atop y\in\mathbb{B}_{\delta}(\bar{y})}}\Pi_{i}(x,y,v) & \subset & \overline{\co}\Pi(\bar{x},\bar{y},v).\label{eq:convexified-target}\end{eqnarray}
\end{lem}
\begin{proof}
First, we note that \eqref{eq:nonconvex-target} implies \eqref{eq:convexified-target}.
If \eqref{eq:nonconvex-target} holds, then by Lemma \ref{lem:nested-sets}
we have \begin{eqnarray*}
\bigcap_{N\in\mathbb{N},\delta>0}\overline{\co}\bigcup_{i>N}\bigcup_{{x\in\mathbb{B}_{\delta}(\bar{x})\atop y\in\mathbb{B}_{\delta}(\bar{y})}}\Pi_{i}(x,y,v) & = & \overline{\co}\bigcap_{N\in\mathbb{N},\delta>0}\cl\bigcup_{i>N}\bigcup_{{x\in\mathbb{B}_{\delta}(\bar{x})\atop y\in\mathbb{B}_{\delta}(\bar{y})}}\Pi_{i}(x,y,v)\\
 & \subset & \overline{\co}\Pi(\bar{x},\bar{y},v).\end{eqnarray*}
Proving \eqref{eq:nonconvex-target} is equivalent to proving the
following: If $u_{i}\in\Pi_{i}(\bar{x}_{i},\bar{y}_{i},v)$ and $u_{i}\to u$,
$\bar{x}_{i}\to\bar{x}$ and $\bar{y}_{i}\to\bar{y}$ as $i\to\infty$,
then $u\in\Pi(\bar{x},\bar{y},v)$. Consider $v\in\mathbb{\ensuremath{R}}^{n}$
and the sequences of functions $\{x_{i}(\cdot)\}_{i=1}^{\infty}$
and $\{p_{i}(\cdot)\}_{i=1}^{\infty}$ constructed from $\{x_{i,j}\}_{{0\leq j\leq i\atop 1\leq i\leq\infty}}$
and $\{p_{i,j}\}_{{0\leq j\leq i\atop 1\leq i\leq\infty}}$ such that
$x_{i}(0)=\bar{x}_{i}$, $x_{i}(T)=\bar{y}_{i}$, $p_{i}(0)=u_{i}$
and $p_{i}(T)=v$. We therefore need to show that $u\in\Pi(\bar{x},\bar{y},v)$. 

For a fixed $t\in[0,T]$, the map \[
(p,\tilde{x},\tilde{y},\tilde{p})\mapsto-D_{x}^{*}F\big(t,x(t)+\tilde{x}\mid x^{\prime}(t)+\tilde{y}\big)(p+\tilde{p})\]
can be checked to be osc (at where $x(t)$ and $x^{\prime}(t)$ are
defined) from the definition of the coderivatives and the fact that
the map $(x,y)\mapsto N_{\scriptsize\gph(F(t,\cdot))}(x,y)$ is osc.
 The map \[
(p,\tilde{x},\tilde{y},\tilde{p})\mapsto-\overline{\co}D_{x}^{*}F\big(t,x(t)+\tilde{x}\mid x^{\prime}(t)+\tilde{y}\big)(p+\tilde{p})\]
is osc since the convex hull operation preserves outer semicontinuity.
(The proof is elementary, and the steps are shown in \cite{chardiff}
for example.)

Suppose $x(\cdot)$ is such that assumption (1) in the statement holds.
Our problem can be solved if we can show that $p(\cdot)$ satisfies
the differential inclusion \[
p^{\prime}(t)\in-\overline{\co}D_{x}^{*}F\big(t,x(t)\mid x^{\prime}(t)\big)\big(p(t)\big).\]
We try to find $\eta_{k}:[0,T]\to[0,\infty)$ such that for all $t\in[0,T]$,
\begin{eqnarray}
 &  & \lim_{k\to\infty}\eta_{k}(t)=0\label{eq:eta-k-t}\\
 & \mbox{ and } & p_{i_{k}}^{\prime}(t)\in-\overline{\co}D_{x}^{*}F\big(t+\eta_{k}(t)\mathbb{B},x(t)+\eta_{k}(t)\mathbb{B}\mid x^{\prime}(t)+\eta_{k}(t)\mathbb{B}\big)\big(p_{i_{k}}(t)+\eta_{k}(t)\mathbb{B}\big).\nonumber \end{eqnarray}
For each $t\in[0,T]$ and $k\in\overline{1,\infty}$, we have $\left\lfloor t/(\Delta t)\right\rfloor (\Delta t)\leq t\leq\left\lfloor t/(\Delta t)+1\right\rfloor (\Delta t)$,
where $\Delta t=T/i_{k}$ and $\left\lfloor \alpha\right\rfloor $
is the greatest integer not more than $\alpha$. For simplicity, we
consider the case where $t/T$ is irrational. From the definitions
of $x_{i_{k}}(\cdot)$ and $p_{i_{k}}(\cdot)$ and \eqref{eq:disc-diff-incl-condn},
we have \[
p_{i_{k}}^{\prime}(t)\in-\overline{\co}D_{x}^{*}F\big(t_{k},x_{i_{k}}(t_{k})\mid x_{i_{k}}^{\prime}(t)\big)\big(p_{i_{k}}(t_{k}+\Delta t)\big),\]
where $t_{k}=\left\lfloor t/(\Delta t)\right\rfloor (\Delta t)$.
To establish the existence of $\eta_{k}(\cdot)$ in \eqref{eq:eta-k-t},
it suffices to show that for each $t$, \[
\max(|t_{k}-t|,\|x_{i_{k}}(t_{k})-x(t)\|,\|x_{i_{k}}^{\prime}(t)-x^{\prime}(t)\|,\|p_{i_{k}}(t_{k}+\Delta t)-p_{i_{k}}(t)\|)\searrow0\mbox{ as }k\nearrow\infty.\]
We first have $x_{i_{k}}^{\prime}(t)\to x^{\prime}(t)$ and $t_{k}\to t$
as $k\to\infty$. Next, since $p_{i_{k}}(\cdot)$ converges uniformly
to $p(\cdot)$, we have \begin{eqnarray}
 &  & \|p_{i_{k}}(t_{k}+\Delta t)-p_{i_{k}}(t)\|\nonumber \\
 & \leq & \underbrace{\|p_{i_{k}}(t_{k}+\Delta t)-p(t_{k}+\Delta t)\|}_{(1)}+\underbrace{\|p(t_{k}+\Delta t)-p(t)\|}_{(2)}+\underbrace{\|p(t)-p_{i_{k}}(t)\|}_{(3)},\label{eq:3-parts-in-p}\end{eqnarray}
so the term on the LHS converges to zero as $k\to\infty$. A similar
argument with $x_{i_{k}}(t_{k})-x(t)$ shows that its norm goes to
zero as $k\to\infty$. So the presence of $\eta_{k}(t)$ satisfying
\eqref{eq:eta-k-t} is established. 

Since $p(\cdot)$ is continuous on the compact set $[0,T]$, it is
uniformly continuous. This implies that for any $\epsilon>0$, we
can find $K$ such that term (2) in \eqref{eq:3-parts-in-p} has norm
less than $\epsilon$ for all $k>K$. The condition that $\eta_{k}(t)\leq\eta(t)$
for all $t\in[0,T]$ for some $\eta(\cdot)\in L_{1}([0,T],\mathbb{R}^{n})$
(in fact, $\eta(\cdot)\in L_{\infty}([0,T],\mathbb{R}^{n})$) follows
easily. All the conditions for Lemma \ref{lem:Smirnov-L4.4} are satisfied,
and we have $u\in\Pi(\bar{x},\bar{y},v)$ as needed.
\end{proof}
Though condition (1) may look more complicated than \eqref{eq:conj-step3}
alone, it can be understood as a measurability condition on $x(\cdot)$
and $p(\cdot)$. We collect the previous results to obtain an estimate
of the convexified limiting coderivative of the reachable map.
\begin{thm}
\label{thm:co-D*F-thm}(Convexified coderivative of reachable map)
The formula \eqref{eq:R-conv-coderv} holds provided:
\begin{itemize}
\item [(a)]For any $\epsilon>0$ and $x\in\mathbb{R}^{n}$, there is some
$I$ such that $\mathbf{d}(R(x),R_{i}(x))<\epsilon$ for all $i>I$.
\item [(b)]For all $x$ and $t$, $F(t,\cdot)$ is locally bounded at $x$.
\item [(c)]There is some $b(\cdot)\in L_{1}([0,T],\mathbb{R}^{n})$ such
that  $\|D_{x}^{*}F(t,x\mid y)\|^{+}\leq b(t)$ for all $(t,x,y)\in[0,T]\times\mathbb{R}^{n}\times\mathbb{R}^{n}$.
\item [(d)]Assumption (1) of Lemma \ref{lem:3rd-containment} holds.
\end{itemize}
\end{thm}
\begin{proof}
This combines Lemma \ref{lem:On-conj-1}, Theorem \ref{thm:discrete-coderv}
and Lemma \ref{lem:3rd-containment}. From (b) and (c), standard methods
of set-valued analysis imply that $F(t,\cdot)$ is locally Lipschitz,
so the requirements for Theorem \ref{thm:discrete-coderv} are satisfied.
The condition $\|D_{x}^{*}F(t,x\mid y)\|^{+}\leq b(t)$ is equivalent
to the condition on $\overline{\co}D_{x}^{*}F(t,x\mid y)$ in Lemma
\ref{lem:3rd-containment}. 
\end{proof}
Conditions (b) and (c) are typical assumptions for (EL), (TC) and
(WP) to hold. Condition (a) is a mild assumption on how the discretized
reachable map can approximate the continuous reachable map, and Condition
(d) relates the discretized paths to continuous paths. The procedure
of passing a sequence of discrete problems to the limit seems to make
it unavoidable that assumption (d) has to hold, and that the conclusion
can only be expressed in terms of convexified limiting coderivatives.
The conditions (EL), (TC) and (WP) are usually proved with direct
methods in analysis rather than through discrete approximations, so
it remains to be seen whether Theorem \ref{thm:co-D*F-thm} can be
further strengthened with such techniques. 
\begin{rem}
(Graph convex $F(t,\cdot)$) The discrete case suggests that when
$F(t,\cdot)$ is graphically convex for all $t$, then \eqref{eq:R-coderv-conj}
is actually an equation. For the continuous case, we study \eqref{eq:R-conv-coderv}
instead, and ask whether \eqref{eq:R-conv-coderv} is an equation
when $F(t,\cdot)$ is graphically convex for all $t$. In this case,
\eqref{eq:conj-step2} is an equation, but equality for \eqref{eq:conj-step1}
requires further assumptions. The reverse inclusion for \eqref{eq:conj-step3}
holds if every continuous path on the RHS can be described as a limit
of sequences on the left hand side. Such results may already be in
the literature. We cite \cite[Theorem 4.16]{Smirnov02} for example,
which states that the reverse inclusion in \eqref{eq:conj-step3}
holds when $F(\cdot,\cdot)$ is independent of its first argument
$t$ and is Lipschitz.
\end{rem}

\section{Conclusion}

In this paper, we study how discrete and differential inclusions depend
on the initial conditions. The advantage of such results over necessary
optimality conditions is that such results give an indication of how
to perturb the initial point to optimality. The results for discrete
inclusions seem quite satisfactory, but the results for differential
inclusions still require further improvement.

\bibliographystyle{alpha}
\bibliography{../refs}

\begin{thebibliography}{Pan11b}

\bibitem[AC84]{AC84}
J.-P. Aubin and A.~Cellina.
\newblock {\em Differential Inclusions: Set-Valued Maps and Viability Theory}.
\newblock Springer, New York, 1984.
\newblock Grundlehren der mathematischen Wissenschaften, Vol 264.

\bibitem[AF90]{AF90}
J.-P. Aubin and H.~Frankowska.
\newblock {\em Set-Valued Analysis}.
\newblock Birkh{\"a}user, Boston, 1990.
\newblock Republished as a Modern Birkh{\"a}user Classic, 2009.

\bibitem[Cla83]{Cla83}
F.H. Clarke.
\newblock {\em Optimization and Nonsmooth Analysis}.
\newblock Wiley, Philadelphia, 1983.
\newblock Republished as a SIAM Classic in Applied Mathematics, 1990.

\bibitem[DL92]{DL92}
A.L. Dontchev and F.~Lempio.
\newblock Difference methods for differential inclusions: A survey.
\newblock {\em SIAM Rev.}, 34(2):263--294, 1992.

\bibitem[DQZ06]{DQZ06}
A.L. Dontchev, M.~Quincampoix, and N.~Zlateva.
\newblock Aubin criterion for metric regularity.
\newblock {\em J. Convex Anal.}, 3:45--63, 2006.

\bibitem[Iof81]{Iof81}
A.D. Ioffe.
\newblock Nonsmooth analysis: differential calculus of non-differentiable
  mappings.
\newblock {\em Trans. Amer. Math. Soc.}, 266:1--56, 1981.

\bibitem[Mor06]{Mor06}
B.S. Mordukhovich.
\newblock {\em Variational Analysis and Generalized Differentiation {I} and
  {II}}.
\newblock Springer, Berlin, 2006.
\newblock Grundlehren der mathematischen Wissenschaften, Vols 330 and 331.

\bibitem[Pan11a]{chardiff}
C.H.J. Pang.
\newblock Characterizing generalized derivatives of set-valued maps: Extending
  the {A}ubin and {M}ordukhovich criterions.
\newblock 2011.
\newblock Available in http://arxiv.org/abs/1106.2338.

\bibitem[Pan11b]{set_diff}
C.H.J. Pang.
\newblock Generalized differentiation with positively homogeneous maps:
  Applications in set-valued analysis and metric regularity.
\newblock {\em Math. Oper. Res.}, 36:3:377--397, 2011.

\bibitem[RW98]{RW98}
R.T. Rockafellar and R.J.-B. Wets.
\newblock {\em Variational Analysis}.
\newblock Springer, Berlin, 1998.
\newblock Grundlehren der mathematischen Wissenschaften, Vol 317.

\bibitem[Smi02]{Smirnov02}
G.~Smirnov.
\newblock {\em Introduction to the theory of differential inclusions}.
\newblock Amer. Math. Soc., Providence, RI, 2002.
\newblock Graduate Studies in Mathematics, Volume 41.

\bibitem[Vin00]{Vinter00}
R.B. Vinter.
\newblock {\em Optimal control}.
\newblock {Birkh{\"a}user}, Boston, 2000.
\newblock Republished as a Modern {Birkh{\"a}user} classic, 2010.

\end{thebibliography}

\end{document}